\documentclass[a4paper,11pt]{article}
\usepackage[latin1]{inputenc}
\usepackage[T1]{fontenc}
\usepackage{a4wide}
\usepackage{fontenc}
\usepackage{color}
\usepackage{amsmath}
\usepackage{graphicx}
\usepackage{array}
\usepackage{amsfonts}
\usepackage{amssymb}
\usepackage{euscript}
\usepackage{delarray}
\usepackage{latexsym}
\usepackage{enumerate}
\usepackage{float}
\usepackage{algorithmicx}
\usepackage{algorithm}
\usepackage{algpseudocode}
\usepackage{url}

\title{Kriging-based sequential design strategies using fast cross-validation techniques with extensions to multi-fidelity computer codes}
\author{Loic Le Gratiet $^\dag$  $^\ddag$ and Claire Cannamela $^\ddag$  \\ \\ $^\dag$ Universit\'e Paris Diderot 75205 Paris Cedex 13 \\ \\ $^\ddag$ CEA, DAM, DIF, F-91297 Arpajon, France  }

\newtheorem{prop}{Proposition}

\begin{document}
\maketitle

\section{Abstract}

Kriging-based surrogate models have become very popular during the last decades  to approximate a computer code output from  few simulations.  In practical applications, it is  very common to sequentially add new simulations to obtain  more accurate approximations. We propose in this paper a  method of kriging-based sequential design which combines both the error  evaluation providing by the kriging model and the observed errors of a Leave-One-Out cross-validation procedure. This method is proposed in two versions, the first one  selects points one at-a-time. The second one  allows us to parallelize the simulations and to add several design points at-a-time. Then, we extend these strategies   to multi-fidelity co-kriging models which allow us to surrogate a complex code using fast approximations of it. The main advantage of these extensions is that it not only provides   the new locations where to perform simulations but also which versions  of code have to be simulated (between the complex one or one of  its fast approximations).

A real multi-fidelity application  is used to illustrate the efficiency of the proposed approaches. In this example,  the accurate code is a two-dimensional finite element model and the less accurate one is a one-dimensional approximation of the system.

\paragraph{Keywords} : kriging, co-kriging, sequential design, cross-validation, resource  allocations, multi-fidelity codes.

\section{Introduction}

Kriging-based surrogate models have become very popular during the last decades to design and analyze computer experiments \cite{SACKS89}. The reader is referred to  the books of  \cite{S99}, \cite{S03} and \cite{R06} for more detail about kriging models. Usually, in real applications, two stages are performed to surrogate a computer code with a kriging model. The first one consists in building a  kriging model from simulations coming from an initial experimental design set. Many methods exist to build the initial design set, in order to ensure appropriate space filling properties,  the reader is referred to \cite{FAN06} for  a non-exhaustive review of them. The second stage   consists in adding simulations sequentially at new design points which complete  the initial one. The selection of the new points are usually based on criteria to improve the global accuracy of the kriging model and this will be our goal in this paper. To be complete, we mention that sequential kriging has also been widely used in optimization (see \cite{DRJ98}, \cite{PiGi12}) and  to estimate  probabilities  of failure \cite{Bec11}

Kriging models are a powerful tool to enrich an experimental design set since it provides through the kriging variance - also called predictor's Mean Squared Error (MSE) or variance of prediction  -  an estimation of the model MSE. Kriging literature provides lot of criteria  usually based on the kriging variance for sequentially design the experiments \cite{SACKS89}. Furthermore, \cite{Bat96} and \cite{Pi10} propose more efficient criteria by considering the Integrated MSE (IMSE). It consists in integrating the mean value of the MSE  integrated over the input parameter space. We note though  that the IMSE  can be computationally expensive to assess, especially when the dimension increases. Although these criteria are efficient  for many  cases, they can suffer from an important flaw when the accuracy of the kriging model is not homogeneous  over the input parameter  space. Indeed, the kriging variance  is determined by  the distances between prediction  and design  points   but not by   the real  model errors. To fix this important flaw, we can use the Empirical IMSE suggested in \cite{SACKS89} which  evaluates the model errors through a test set. Nevertheless, in a complex computer code  framework, it could be too   expensive to consider an external test set and  cross-validation (CV) based criteria are more significant. As an illustration \cite{Klei04} and \cite{Klei06}   combine a bootstrapping and a CV procedure to evaluate the predictor's MSE. Although this method improves  the classical approach, it still does not take into account the real model errors. We note that a strength of the method proposed by  \cite{Klei04}  is that it can be applied  to others type of surrogate models than the kriging one.

The first focus of this paper is  on sequential design to improve the accuracy of a kriging model. In particular, we propose new criteria   combining  the kriging variance and the Leave-One-Out CV (LOO-CV) errors. The CV errors   allow for focusing  the new observations on regions where the real model errors are large. Furthermore, thanks to   the equations presented in \cite{Dub83}, the LOO-CV equations are fast to compute and thus the suggested approach is not expensive.

Kriging has recently  been extended to allow for the use of coarse versions of a complex computer code to improve the accuracy of its approximation. These  so-called multi-fidelity surrogate models  have  become of growing interest. A first one was proposed in \cite{Crai98} which is based on a linear regression formulation and was improved by \cite{CG09}  through a Bayes linear formulation. The first multi-fidelity model using   co-kriging   was suggested by \cite{KO00}. Then, several works dealing with this model have been developed \cite{KO01}, \cite{Hig04}, \cite{Ree04}, \cite{QW07}. Defining   sequential design strategies  in a multi-fidelity framework  is also of interest and is still an open problem. A method based on nested Latin hypercube designs is suggested in  \cite{Xio12}. However, it does not allow  for adding few simulations (e.g. it cannot perform an one step at-a-time sequential design)  and it does not take into account the accuracy of the coarse code versions and the time ratios between two code levels. 

 The second focus of this paper is on  sequential design for co-kriging model. We adapt the new strategies suggested for the kriging model to the multi-fidelity co-kriging one. The strength of the proposed extensions is that they not only provide  the  new  points  where to perform new simulations but  they also determine  which version of code is worth being simulated. These new criteria take into account  the computational time ratios between code versions. They  are based on  a proxy of the  IMSE reduction and on an original result giving the contribution  of each code on the total variance of the model. We note that sequential design in a multi-fidelity framework has also been applied for optimization purposes \cite{AF07} and \cite{Hu06}.

The paper is organized as follows. First, we introduce the kriging model  and   present our CV-based sequential design strategies. We illustrate these   strategies in   tabulated  functions. Secondly, we present the proposed co-kriging multi-fidelity model and the extensions of the previous strategies. Finally, we apply  the sequential co-kriging approach to a mechanical example. 

\section{Kriging models and sequential designs}\label{section1}

In this Section, we briefly introduce the kriging model and some of its classical sequential design criteria. Then, we will present our   sequential strategies to enhance kriging models considering the region with  large LOO-CV errors.

\subsection{The Kriging model}\label{subsection11}

The kriging model is a widely used method to surrogate the output of a computer code from few simulations  \cite{SACKS89}. Let us denote by $y(x)$ the output of the code at point  $x \in Q$. Here, $y(x)$ is a scalar and   $Q \subset \mathbb{R}^d$. Furthermore, we denote by $\mathbf{D} = \{x_1,\dots,x_n\}$ the experimental design set and $\mathbf{y}_n = y(\mathbf{D})$ the value of $y(x)$ at points in $\mathbf{D}$.

In a kriging framework, we set that the prior knowledges about the code can be modeled by a Gaussian process $Y_0(x)$. Usually, we consider a Gaussian process  with mean of the form $m_0(x) = \mathbf{f}'(x)\beta$, with $\mathbf{f}'(x) = (f_1(x),\dots,f_p(x))$ and with covariance function $k_0(x,\tilde{x}) = \sigma^2r\left( x,\tilde{x} ;\theta \right)$. Then, the kriging equations are given by the Gaussian process $Y_0(x)$ conditioned by its known values $\mathbf{y}_n$ at points in $\mathbf{D}$:
\begin{equation}\label{eq1}
Y_n(x) = [Y_0(x) | Y_0(\mathbf{D}) = \mathbf{y}_n] \sim  \mathcal{GP} \left(m_n(x) , k_n(x,\tilde{x})\right)
\end{equation}
where:
\begin{equation}\label{eq2}
m_n(x) =  \mathbf{f}'(x) \hat{\beta} + \mathbf{r}' (x)\mathbf{R}^{-1}(\mathbf{y}_n - \mathbf{F}\hat{\beta} )
\end{equation}
and:
\begin{equation}\label{eq3}
k_n(x,\tilde{x}) = \sigma^2 \left( r(x,\tilde{x})  -  \left( \begin{array}{cc} \mathbf{f}'(x) &  \mathbf{r}' (x)  \end{array} \right)   \left( \begin{array}{cc} 0 & \mathbf{F}' \\ \mathbf{F}   & \mathbf{R}  \end{array} \right) ^{-1} \left( \begin{array}{c} \mathbf{f}(x) \\ \mathbf{r} (x) \end{array} \right) \right)
\end{equation}
where $'$ stands for the transpose, $ \mathbf{F}  $ are the values of $\mathbf{f}'(x)$ at points in $\mathbf{D}$, $\mathbf{r}(x)$ is the correlation vector between $\mathbf{D}$ and $x$ with  respect to the correlation function $r(x,\tilde{x})$, $\mathbf{R}$ is the correlation  matrix of $\mathbf{D}$  with respect to  $r(x,\tilde{x})$ and $ \hat{\beta} = ( \mathbf{F}' \mathbf{ R}^{-1}\mathbf{ F})^{-1}\mathbf{F}' \mathbf{R}^{-1}\mathbf{y}_n  $  is the usual least-squares estimates of $\beta$. The model parameters $\sigma^2$ and $\theta$ could be estimated by maximizing their Likelihood \cite{S03} or with a cross-validation procedure \cite{R06}. Furthermore, the Maximum Likelihood Estimate (MLE) of $\sigma^2$ is given by $\hat{\sigma}^2 = (\mathbf{y}_n-\mathbf{F} \hat{\beta})'\mathbf{R}^{-1}(\mathbf{y}_n-\mathbf{F} \hat{\beta})/(n-p) $.

\paragraph{ 1 point  at-a-time Sequential design\\}Now, let us suppose that we want to add a new point $x_{n+1}$  in $\mathbf{D}$ in order to enhance the accuracy of the kriging model.
From the  kriging variance $k_n(x,x)$ - representing the model MSE -  some sequential design methods have been derived \cite{SACKS89}, \cite{Bat96} and   \cite{Pi10}.  
A first one consists in adding $x_{n+1}$ where the kriging variance is the largest (see \cite{SACKS89}):
\begin{equation}\label{eq3bis}
x_{n+1} = \arg \max_x  k_n(x,x) 
\end{equation}
However, as presented in \cite{Klei04},  its performance is   poor. Then, it has been improved with a criterion which  consists in adding the new point which leads the most important IMSE reduction (see \cite{Bat96} and   \cite{Pi10}):
\begin{equation}\label{eq4}
x_{n+1} = \arg \max_x  \int_{u \in Q} k_n(u,u) -  k_{n+1}(u,u) \, du
\end{equation}
Here, the covariance kernels $k_{n+1}(u,\tilde{u})$   corresponds to the one of the Gaussian process $Y_n(u)$ (\ref{eq1}) conditioned by a new observation at $x$.
Furthermore, the equation (\ref{eq3}) shows that the kriging variance does not depend on the observations if we consider  known  the  parameters $\sigma^2$ and $\theta$.
Therefore, in that case, $k_{n+1}(u,u)$  can be computed without new simulations. We denote by MinIMSE this criterion. Finally, we also consider the  criterion presented by \cite{Klei04} using a Jackknife estimator for the predictor's variance. Its principle is the following one. Let us consider $m_{n,-i}(x)$ the kriging mean built without the i$^{th}$ observation, the jackknife variance is given by:
\begin{equation}
s^2_{jack}(x) = \frac{1}{n(n-1)} \sum_{i=1}^{n} (\tilde{y}_i - \bar{\tilde{y}})^2
\end{equation} 
where $\tilde{y}_i =  nm_n(x) - (n-1) m_{n,-i}(x)$ and $\bar{\tilde{y}} = \sumç_{i=1}^n  \tilde{y}_i / n$. Then, we consider candidate points coming from   a maximin LHS Design \cite{FAN06} and we add those which maximize the jackknife variance. We denote by KleiCrit this criterion.

\paragraph{$\mathbf{q}$ points at-a-time Sequential design \\}
There is a natural way to extend these algorithms when the simulations can be performed simultaneously. Indeed, the covariance kernel $k_{n+1}(x,\tilde{x})$ of the Gaussian process $Y_n(x)$  conditioned by   the new observation at   point $x_{n+1}$  can be computed without knowing $y(x_{n+1})$ when we consider the model parameters $\sigma^2$ and $\theta$ as known. Then, from $k_{n+1}(x,\tilde{x})$, we can find a new point $x_{n+2}$ where to perform a new simulation using the same criterion as in equation (\ref{eq4}) and the kernel $k_{n+2}(x,\tilde{x})$. Thus, considering the parameters   $\sigma^2$ and $\theta$ as known (they are fixed to their estimated values), we can determine with this procedure $q$ good locations where to perform simulations. We call this method  the  liar sequential kriging.  This idea is also presented in a framework of kriging-based optimization in \cite{Ginsb10}.

\subsection{LOO-CV based  strategies for kriging sequential design}\label{subsection12}

We present in this subsection  new  sequential-kriging strategies. The main difference between these new strategies and the previous ones is that they  take into account the real model errors through the LOO-CV equations. 

The proposed    sequential methods are based on the following original proposition. It gives  the closed form expressions of the LOO-CV equations where  the model  parameters $\beta$ and $\sigma^2$ are re-estimated after each removed point. This result is already known when $\sigma^2$ is fixed \cite{Dub83} and \cite{R06}.

\paragraph{Notations:} $  \mathbf{A}_{i,i}$ is the  i$^{th}$ element of the main diagonal of $ \mathbf{A}$, $  \mathbf{A}_{i }$ is the  i$^{th}$ row of the matrix  $\mathbf{A}$, $  \mathbf{A}_{-i }$ is the matrix $\mathbf{A}$ without its   i$^{th}$ row,  $  \mathbf{A}_{-i,i}$ is the  i$^{th}$ column of $ \mathbf{A}$ without its i$^{th}$ element, $ \mathbf{A}_{i,-i} =   \mathbf{A}_{-i,i}'$ and $  \mathbf{A}_{-i,-i}$  is the matrix  $ \mathbf{A}$ without the  i$^{th}$ row and column.

\begin{prop}\label{prop1}
Let us denote by $Y_{n,-i}(x)$ the Gaussian process $Y_0(x)$ conditioned by the values $\mathbf{y}_{n,-i} = y(\mathbf{D})  \setminus y(x_i) $.  Then, the predictive mean   of $Y_{n,-i}(x)$ at point $x_i$ is given by:
\begin{equation}\label{eq5}
m_{n,-i}(x_i) = y(x_i)  -  \left[ \mathbf{R}^{-1}( \mathbf{y}_{n,-i} - \mathbf{F}_{-i}\hat{\beta}_{-i} ) \right]_i / \left[ \mathbf{R}^{-1}\right]_{i,i}
\end{equation}
where
$\hat{\beta}_{-i} = (\mathbf{F}'_{-i} \mathbf{K}_i   \mathbf{F}_{-i} )^{-1}\mathbf{F}'_{-i} \mathbf{K}_i \mathbf{y}_{n,-i}$
and $\mathbf{K}_i = \left[ \mathbf{R}^{-1}\right]_{-i,-i}  - \left[ \mathbf{R}^{-1}\right]_{-i, i}\left[ \mathbf{R}^{-1}\right]_{ i,-i} / \left[ \mathbf{R}^{-1}\right]_{i,i}$.
Furthermore, the predictive variance of $Y_{n,-i}(x)$ at point $x_i$ is given by:
\begin{equation}\label{eq6}
k_{n,-i}(x_i) =  \sigma_{-i}^2 / \left[ \mathbf{R}^{-1}\right]_{i,i} + \varsigma_{-i}(x_i)
\end{equation}
where
$
\varsigma_{-i}(x_i) = \left(\left[ \mathbf{R}^{-1} \mathbf{F} \right]_i/ \left[ \mathbf{R}^{-1}\right]_{i,i} \right)'(\mathbf{F}'_{-i} \mathbf{K}_i   \mathbf{F}_{-i} )^{-1}\left(\left[ \mathbf{R}^{-1} \mathbf{F} \right]_i/ \left[ \mathbf{R}^{-1}\right]_{i,i} \right)
$ 
and \\
$
  \sigma_{-i}^2 = \left( \mathbf{y}_{n,-i} - \mathbf{F}_{-i}\hat{\beta}_{-i} \right)' \mathbf{K}_i \left( \mathbf{y}_{n,-i} - \mathbf{F}_{-i}\hat{\beta}_{-i} \right)/(n-p-1)
$
\end{prop}

The previous proposition provides a powerful tool to compute the LOO-CV predictive means and variances. Indeed, several elements of the equations have been already computed during the models construction (e.g. the inverse of the matrix $\mathbf{R}$). Consequently, the LOO-CV equations are fast to compute and can be easily recomputed at each step of the sequential strategy. We note that the originality of this result is the estimation of $ \sigma_{-i}^2$. As we use the value of $k_{n,-i}(x_i)$ strongly depending on $ \sigma_{-i}^2$  in our forthcoming developments,  it is important to well estimate it.

Now, let us denote by $\mathbf{e}^2_{\mathrm{LOO-CV}} = \left[( \left(  y(x_i) - m_{n,-i}(x_i)  \right)^2\right]_{i=1,\dots,n}$ the vector of the LOO-CV squared errors and $\mathbf{s}^2_{\mathrm{LOO-CV}} =\left[ k_{n,-i}(x_i)  \right]_{i=1,\dots,n}$ the vector of the LOO-CV variances. Furthermore, let us consider the Voronoi cells $(V_i)_{i=1,\dots,n}$ associated with the points $(x_i)_{i=1,\dots,n}$:
\begin{equation}\label{eq7}
V_i = \{ x \in Q , \, ||x-x_i|| \leq ||x-x_j||, \, \forall j \neq i \}, \, i,j = 1,\dots,n
\end{equation}
In the remainder  of this section, we present two strategies to sequentially add simulations which use $\mathbf{e}^2_{\mathrm{LOO-CV}}$, $\mathbf{s}^2_{\mathrm{LOO-CV}}$ and $V_i$. The intuitive idea of the suggested criteria is to enhance the predictive variance in the locations where the LOO-CV errors are important. 

\paragraph{LOO-CV-based 1 point  at-a-time Sequential design\\}
Let us denote by $x_{n+1}$ the new point  that we want to add to $\mathbf{D}$. We consider the point solving the following problem:
\begin{equation}\label{eq8}
x_{n+1} = \arg \max_x  k_n(x,x)\left(1+\sum_{i=1}^n { \frac{[\mathbf{e}^2_{\mathrm{LOO-CV}}]_i}{[\mathbf{s}^2_{\mathrm{LOO-CV}}]_i} \mathbb{I}_{x\in V_i} } \right)
\end{equation}
This criterion considers the predictor's MSE $k_n(x,x)$ adjusted with the LOO-CV errors and variances. For equivalent $k_n(x,x)$, the criterion favors the points close to an experimental design point with  large LOO-CV errors.  Furthermore, if two  points are  in the same Voronoi cell, the one with the largest predictor's MSE is considered. Therefore, a sequential strategy with this criterion focus on the regions of $Q$ where the LOO-CV errors are the largest. We note that the standardization with $\mathbf{s}^2_{\mathrm{LOO-CV}}$ is important since it is not necessary to enlarge the predictor's MSE in the regions where it is well or over estimated.

\paragraph{LOO-CV-based  $\mathbf{q}$  points  at-a-time Sequential design\\}
We extend here the previous criterion for a $q$ points at-a-time sequential design. First, we emphasize that the liar sequential kriging is not relevant for this new criterion. Indeed, conditioning on model parameters,  with a liar  method we can compute the kriging variances $(k_{n+i}(x,x))_{i=1,\dots,q}$ but not  the LOO-CV equations.  Therefore, we use another strategy to propose $q$ new locations where to perform the simulations. This approach is proposed in \cite{Dub11} in a different framework. The idea of the suggested method is to select the $q$ best points with respect to the criterion (\ref{eq8}) from $N$ candidate points. These $N$ candidate points are  chosen with the following algorithm. 
\begin{enumerate}
\item Generate $N_{\mathrm{MCMC}}$ samples  with respect to the  probability density function  proportional to $k_n(x,x)$ with a suitable Markov Chain Monte Carlo (MCMC) technique  \cite{RobCas04}.
\item Extract from these samples $N$ representative points with a $N$-means clustering technique \cite{Mac67}.
\end{enumerate}

As presented in  \cite{Dub11}  the use of this algorithm to select $N$ candidate points in a kriging framework is   efficient. Indeed, it allows us to concentrate the points at the modes of the kriging variance. In the proposed strategy,  we always take $N \geq q$ and we choose from the $N$ cluster centers $(C_i)_{i=1\dots,N}$ the $q$ points where $k_{n,\mathrm{adj}}(x,x) =  k_n(x,x)\left(1+\sum_{i=1}^n { \frac{[\mathbf{e}^2_{\mathrm{LOO-CV}}]_i}{[\mathbf{s}^2_{\mathrm{LOO-CV}}]_i} \mathbb{I}_{x\in V_i} } \right)$ is the largest. For the MCMC procedure, we use a Metropolis-Hastings (M-H) algorithm with a Gaussian jumping distribution.  It is centered on the last sample point and has a standard deviation such that the acceptance rate is around 30\% (see \cite{RobCas04}). Furthermore, we set $N_\mathrm{MCMC}$ such that $N_\mathrm{MCMC} \gg N$. For the $N$-means procedure, we choose the value of $N$ with respect to the following criterion:
\begin{equation}\label{eq9}
\max_{N\geq q } \min_{x \in (C_i)_i}  k_n(x,x)
\end{equation}
where $(C_i)_{i=1\dots,N}$ are the cluster centers. This criterion prevents from having a cluster center  in a region where the kriging variance is close to zero. Furthermore, if the number of clusters is too high, the cluster centers get away from the modes and consequently the value of $\min_{x \in (C_i)_{i=1\dots,N}}  k_n(x,x)$ decreases. Therefore, this criterion also prevents from having a number of clusters too large.

\section{Sequential design in a multi-fidelity framework}\label{section2}

A computer code can often be run at different levels of accuracy. In this case, it can be worth using low-accuracy versions of  a code to improve   its approximation. Co-kriging models are well suited to build such multi-fidelity surrogate models  (see \cite{KO00} and  \cite{QW07}). In this section, we present the considered multi-fidelity co-kriging models and we  extend the previous  sequential design strategies in this framework.  We note that, in a multi-fidelity framework,  the search for  the best locations where to run the code is not the only point of interest. Indeed, once the best  locations are determined, we also have to decide which  code level   is worth being run. This will not only depend  on the time-ratios between the code levels but also on the contribution  of each code level to the total predictor's MSE.

\subsection{Multi-fidelity co-kriging models}\label{subsection21}

Let us suppose that we want to surrogate a  computer code output $y^s(x)$   and that coarse versions of this code $(y^l(x))_{l=1,\dots,s-1}$ are available. These codes are sorted by order of fidelity from the less accurate $y^1(x)$ to the most  accurate $y^{s-1}(x)$. All code levels are modeled by Gaussian processes $(Y^l(x))_{l=1,\dots,s}$ with respect to the following relationship with $l=2,\dots,s$:
\begin{equation}\label{eq10}
\left\{
\begin{array}{l}
Y^l(x) = \rho_{l-1}\tilde{Y}^{l-1}(x)+\delta^l(x) \\
\tilde{Y}^{l-1}(x) \perp \delta^l(x) \\
\tilde{Y}^{l-1}(x)  \sim [Y^{l-1}(x)| \mathbf{Y}^{(l-1)} = \mathbf{y}^{(l-1)}]
\end{array}
\right.
\end{equation}
where $\mathbf{Y}^{(l)} = (\mathbf{Y}^1, \dots, \mathbf{Y}^l)$ with   $\mathbf{Y}^l = Y^l(\mathbf{D}^l)$, $\mathbf{y}^{(l)} = (\mathbf{y}^1, \dots, \mathbf{y}^l)$ with $\mathbf{y}^l = y^l(\mathbf{D}^l)$ and $(\mathbf{D}^l)_{l=1,\dots,s}$ are the experimental design sets at level $l$ with $n_l$ points and such that  $\mathbf{D}^s \subseteq \mathbf{D}^{s-1}  \subseteq \dots \subseteq \mathbf{D}^1$. We note that the nested property is not necessary to build the model but allows for  efficient parameter estimations. Furthermore, conditioning on parameters $\beta_l$, $\sigma^2_l$ and $\theta_l$,  $\delta^l(x)$ is a Gaussian process of mean $\mathbf{f}_ l'(x)\beta_l$ with $\mathbf{f}_l'(x) = (f_1^l(x),\dots,f_{p_l}^l(x))$ and  covariance kernel $\sigma_l^2r_l(x,\tilde{x},\theta_l)$.  By convention $Y^1(x)$ has the same distribution as $\delta^1(x)$. 

The multi-fidelity co-kriging model at level $l = 2,\dots,s$ is given by the following distribution:
\begin{equation}\label{eq11}
Y^l_{n_l}(x) = [Y^l(x)|\mathbf{Y}^{(l)} = \mathbf{y}^{(l)} ] \sim \mathcal{GP}(\mu_{n_l}^l(x),k_{n_l}^l(x,\tilde{x}))
\end{equation}
with:
\begin{equation}\label{eq12}
\mu_{n_l}^l(x) = \hat{\rho}_{l-1}\mu_{n_{l-1}}^{l-1}(x) +\mathbf{f}_l(x)\hat{\beta}_l  +  \mathbf{r}_l'(x)\mathbf{R}_l^{-1}(\mathbf{y}^l-\mathbf{F}_l\hat{\beta}_l - \hat{\rho}_{l-1}y^{l-1}(\mathbf{D}^{l}))
\end{equation}
and:
\begin{equation}\label{eq13}
k_{n_l}^l(x,\tilde{x}) = \hat{\rho}_{l-1}^2k_{n_{l-1}}^{l-1}(x,\tilde{x}) +  \sigma_l^2 \left( r_l(x,\tilde{x})  -  \left( \begin{array}{cc} \mathbf{h}_l'(x) &  \mathbf{r}_l' (x)  \end{array} \right)   \left( \begin{array}{cc} 0 & \mathbf{H}_l' \\ \mathbf{H}_l   & \mathbf{R}_l  \end{array} \right) ^{-1} \left( \begin{array}{c} \mathbf{h}_l(\tilde{x}) \\ \mathbf{r}_l (\tilde{x}) \end{array} \right) \right)
\end{equation}
where $ \mathbf{F}_l  $ are the values of $\mathbf{f}'_l(x)$ at points in $\mathbf{D}^l$, $\mathbf{r}_l(x)$ is the correlation vector between $\mathbf{D}^l$ and $x$ with  respect to the correlation function $r_l(x,\tilde{x})$, $\mathbf{R}_l$ is the correlation  matrix of $\mathbf{D}^l$  with respect to  $r_l(x,\tilde{x})$, $\mathbf{h}_l'(x) = \left[\mu_{n_{l-1}}^{l-1}(x)  \, \,  \mathbf{f}'_l(x) \right] $  and $\left( \begin{array}{c} \hat{\rho}_{l-1} \\  \hat{\beta}_l \end{array} \right) = ( \mathbf{H}_l' \mathbf{ R}_l^{-1}\mathbf{H}_l)^{-1}\mathbf{H}_l' \mathbf{R}_l^{-1}\mathbf{y}^l  $  with $\mathbf{H}_l = [y^{l-1}(\mathbf{D}^l) \, F_l]$ is the usual least-squares estimates of $\left( \begin{array}{c} \rho_{l-1} \\ \beta_l \end{array} \right)$. Furthermore, the restricted maximum likelihood estimate of $\sigma_l^2$ is given by $\hat{\sigma}_l^2 = \left(\mathbf{y}^l - \mathbf{H}_l\left( \begin{array}{c} \hat{\rho}_{l-1} \\  \hat{\beta}_l \end{array} \right)\right)' \mathbf{R}_l^{-1}  \left(\mathbf{y}^l - \mathbf{H}_l\left( \begin{array}{c} \hat{\rho}_{l-1} \\  \hat{\beta}_l \end{array} \right)\right) / (n_l-p_l-1)$ \cite{S03}.

\paragraph{Remark.}
The important property of this co-kriging model is that its MSE  (\ref{eq13}) provides through the term $\rho_{l-1}^2k_{n_{l-1}}^{l-1}$  the contribution of the code level $l-1$ to the total predictor's MSE. Therefore, it can  allow us to determine which  code level  is worth being simulated at a new location $x$.

\subsection{Sequential design for multi-fidelity co-kriging models}\label{subsection23}

The aim of this subsection is to extend  the sequential kriging strategies proposed in Subsection \ref{subsection12} to the suggested  multi-fidelity co-kriging model. These extensions are based on the variance decomposition property presented in Subsection \ref{subsection21} in equation (\ref{eq13}) and on the following  extension of Proposition \ref{prop1}:

\begin{prop}\label{prop5}
For $l=1,\dots,s$, let us consider the multi-fidelity model  (\ref{eq10}) with  $\mathbf{D}^l  \subseteq \mathbf{D}^{l-1}  \subseteq \dots \subseteq \mathbf{D}^1$. We denote by $\xi_j$ the index of $\mathbf{D}^j$  corresponding to the i$^{th}$ point  $x_i^l$ of $\mathbf{D}^l$ with $1 \leq j \leq l $. Then, if we note $\varepsilon_{\mathrm{LOO-CV},l} (x_i^l)$  the LOO-CV error  (i.e. real value minus predicted value) at level $l$  and  point $x_i^l$, we have:
\begin{equation}\label{eq18}
\begin{array}{lll}
\varepsilon_{\mathrm{LOO-CV},l} (x_i^l)  & = &  \hat{\rho}_{l-1}  \varepsilon_{\mathrm{LOO-CV},l-1}  (x_i^l)  +
\left[ \mathbf{R}_l^{-1}\left(
\mathbf{y}^l - \mathbf{H}_l \left( \begin{array}{c} \hat{\rho}_{l-1} \\  \hat{\beta}_l \end{array} \right)
\right) \right]_{\xi_l } \\
\end{array}
/\left[ \mathbf{R}_l^{-1} \right]_{\xi_l,\xi_l}
\end{equation}
where $\varepsilon_{\mathrm{LOO-CV},1}  (x_i^l) $ is given by Proposition \ref{prop1},
$\left( \begin{array}{c} \hat{\rho}_{l-1} \\  \hat{\beta}_l \end{array} \right) = ( \mathbf{H}_{l,-\xi_l}' \mathbf{ K}_l \mathbf{H}_{l,-\xi_l})^{-1}\mathbf{H}_{l,-\xi_l}' \mathbf{K}_l \mathbf{y}^l_{-\xi_l}  $
 and 
$
\mathbf{ K}_l = \left[ \mathbf{R}_l^{-1} \right]_{-\xi_l,-\xi_l}-\left[ \mathbf{R}_l^{-1} \right]_{-\xi_l,\xi_l}\left[ \mathbf{R}_l^{-1} \right]_{\xi_l,-\xi_l} / \left[ \mathbf{R}_l^{-1} \right]_{\xi_l,\xi_l}
$.

Furthermore, if we note $\sigma_{\mathrm{LOO-CV},l}^2 (x_i^l) $ the variance of the   LOO-CV, we have:
\begin{equation}\label{eq19}
\sigma_{\mathrm{LOO-CV},l}^2 (x_i^l) = \hat{\rho}_{l-1}^2 \sigma_{\mathrm{LOO-CV},l-1}^2 (x_i^l) + \sigma_{l,-\xi_l}^2 / \left[ \mathbf{R}_l^{-1} \right]_{\xi_l,\xi_l}
+\varsigma_l
\end{equation}
with $
\varsigma_l =
\mathrm{u}_l^2
\left( \mathbf{H}_{l , -\xi_l }' \mathbf{K}_l \mathbf{H}_{l , -\xi_l } \right)^{-1}
$,  
 $\mathrm{u}_l  =  
\left[ \mathbf{R}_l^{-1}\mathbf{H}_l \right]_{\xi_s } /  [\mathbf{R}_l^{-1}]_{\xi_l,\xi_l}$ and:
\begin{equation}\label{eq20}
 \sigma_{l,-\xi_l}^2 = \frac{\left(\mathbf{y}^l_{-\xi_l} - \mathbf{H}_{l , -\xi_l }\left( \begin{array}{c} \hat{\rho}_{l-1} \\  \hat{\beta}_l \end{array} \right) 
  \right)' \mathbf{K}_s \left(\mathbf{y}^l_{-\xi_l} -  \mathbf{H}_{l , -\xi_l }
\left( \begin{array}{c} \hat{\rho}_{l-1} \\  \hat{\beta}_l \end{array} \right)
 \right)  }{n_l-p_l-2}
\end{equation}
where  $\mathbf{H}_{l , -\xi_l } =[\mathbf{y}^l_{-\xi_l}  \quad \mathbf{F}_{l , -\xi_l }] $.

\end{prop}

Proposition \ref{prop5} provides closed form expressions for the LOO-CV errors and variances. 
From them, the LOO-CV equations are fast to compute and consequently they can be used in a sequential procedure with a low computational cost.
 Furthermore, since the experimental design sets are nested, we state that during the LOO-CV procedure at level $l$, the points are removed from all code levels. Finally, from these  equations, we can adjust the co-kriging variances $\left( k^l_{n_l}(x,\tilde{x}) \right)_{l=1,\dots,s}$ at each level using  the  same method as presented in equation (\ref{eq8}).

\paragraph{1 point  at-a-time sequential co-kriging.}

First, let us consider $x_{\mathrm{new}}$ the point solving the problem:
\begin{equation}\label{eq21}
x_{\mathrm{new}} = \arg \max_x  k_{n_s}^s(x,x)
\end{equation}
Therefore, we want to compute a new simulation at point where the predictor's MSE is maximal.  Now, let us consider two successive code levels $l-1$ and $l$. The question of interest is to estimate which of these two code  levels  is worth being simulated. 

First, thanks to the equation (\ref{eq13}), we can deduce the contribution of each code levels to the predictor's MSE. Let us define the following notation for $l = 2,\dots,s$:
\begin{equation}\label{eq22}
\sigma_{\delta^l}^2(x) = \sigma_l^2 \left( 1   -  \left( \begin{array}{cc} \mathbf{h}_l'(x) &  \mathbf{r}_l' (x)  \end{array} \right)   \left( \begin{array}{cc} 0 & \mathbf{H}_l' \\ \mathbf{H}_l   & \mathbf{R}_l  \end{array} \right) ^{-1} \left( \begin{array}{c} \mathbf{h}_l(x) \\ \mathbf{r}_l (x) \end{array} \right) \right)
\end{equation}
and  $\sigma_{\delta^1}^2(x) = k^1_{n_1}(x,x)$. Then, we have:
\begin{equation}\label{eq23}
k^l_{n_l}(x,x)  =  \sum_{i=1}^l \sigma_{\delta^i}^2(x) \prod_{j=i}^{l-1}\rho_j^2
\end{equation}
Let us consider that the parameters $(\theta_l)_{l=1,\dots,s}$ define  the characteristic length-scales of the kernels $((r_l(x,  \tilde{x};\theta_l))_{i=1,\dots,s}$  (see \cite{R06} p.83). Then, we can approximate the reduction of the IMSE  after adding a  new point $x_\mathrm{new}$ at level $l$  with the following formula:
\begin{equation}\label{eq24}
 \mathrm{IMSE}_\mathrm{red}^l(x_\mathrm{new}) =  \sum_{i=1}^l \sigma_{\delta^i}^2(x_\mathrm{new}) \prod_{j=i}^{l-1}\rho_j^2 \prod_{m=1}^d \theta_i^m
\end{equation}
with $\theta_l = (\theta_l^1,\dots,\theta_l^d)$. Indeed, at each stage,  $\sigma_{\delta^i}^2(x_\mathrm{new})\prod_{j=i}^{l-1}\rho_j^2$ represents 
the contribution of  the bias $\delta^i(x)$ to the co-kriging variance and $\prod_{m=1}^d \theta_i^m$ represents  the volume  of influence of  $x_\mathrm{new}$ at  level $j$.
This criterion is justify by the fact that the reduction of $\mathrm{IMSE}^l$ defined by  $\mathrm{IMSE}^l = \int_Q \sigma_{\delta^l}^2(x)  \,dx $  after adding a new point $x_\mathrm{new}$ has the same order of magnitude than $\sigma_{\delta^i}^2(x_\mathrm{new})$ times the volume of influence $ \prod_{m=1}^d \theta_i^m$ of $x_\mathrm{new}$.

Now, let us consider  that the ratio of computational times  between the codes $y^l(x)$ and $y^{l-1}(x)$ equals $B_{l/l-1}$. It is worth running the code $y^{l-1}(x)$ if  $ B_{l/l-1} \mathrm{IMSE}_\mathrm{red}^{l-1}(x_\mathrm{new})  >  \mathrm{IMSE}_\mathrm{red}^l(x_\mathrm{new})  $, i.e. if the potential uncertainty reduction by running $B_{l/l-1}$ times $y^{l-1}(x)$ is greater than the one when we run one simulation on $y^l(x)$. From this criterion, we can deduce the following  algorithm for an one at-a-time sequential co-kriging model taking into account both the computational ratios between the different code levels and the contribution of each level  to  the total co-kriging variance.

\begin{algorithm}[H]
\caption{One point  at-a-time sequential co-kriging}
\label{algo1}
\begin{algorithmic}[1]
\State Find $x_\mathrm{new}$ such that $x_{\mathrm{new}} = \arg \max_x  k_{n_s}^s(x,x)$
\For { $l =2,\dots,s$}
\If { $\left( \sigma_{\delta^l}^2(x_{\mathrm{new}})  < \mathrm{IMSE}^l \right)$}
\State Run $y^{l-1}(x_{\mathrm{new}})$
\State \textbf{end for}
\Else
\If {$\left( \mathrm{IMSE}_\mathrm{red}^{l-1}(x_\mathrm{new}) / \mathrm{IMSE}_\mathrm{red}^{l}(x_\mathrm{new}) > 1 / B_{l/l-1} \right)$}
\State Run $y^{l-1}(x_{\mathrm{new}})$
\State \textbf{end for}
\EndIf
\EndIf
\EndFor
\If { $ \left( l = s \right)$}
\State Run $y^{l}(x_{\mathrm{new}})$
\EndIf
\end{algorithmic}
\end{algorithm}

\paragraph{Remarks:}   Algorithm  \ref{algo1} evaluates for two successive code levels $l-1$ and $l$, which one is worth being simulated. It starts with the levels one and two, then two and three and so on. When it finds that the level $l-1$ is more promising  than the level $l$, it stops the loop and simulate $x_\mathrm{new}$ at code levels $y^1(x),\dots,y^{l-1}(x)$. Since the loop is defined from level 1 to level $s$, it favors simulations at low code levels. Therefore, it will tend to learn the coarse code versions before learning the accurate ones.  We note that during the loop of the algorithm \ref{algo1}, the parameters are not re-estimated. In fact, they are re-estimated after adding the new point $x_\mathrm{new}$. Moreover, the first test $\sigma_{\delta^l}^2(x_{\mathrm{new}})  < \mathrm{IMSE}^l $ checks in averaged if the code level $l$  at point $x_\mathrm{new}$ is worth being run. Then, the test $\mathrm{IMSE}_\mathrm{red}^{l-1}(x_\mathrm{new}) / \mathrm{IMSE}_\mathrm{red}^{l}(x_\mathrm{new}) > 1 / B_{l/l-1}$ evaluates which code levels between  $l$ and $l-1$ is the most promising. Finally, if we consider that the code level $l$ is more promising than the code level $l-1$, we confront it to the following code level $l+1$. We note that the algorithm  \ref{algo1} is reiterated until   a prescribed accuracy is reached  or the computational time budget is spent.

\paragraph{1 point at-a-time sequential co-kriging with adjusted predictor's MSE.}

From Proposition \ref{prop5},  Algorithm \ref{algo1} and Equation (\ref{eq24}), we can extend the criterion (\ref{eq8}) to the multi-fidelity co-kriging model. Let us consider the following quantity:
\begin{equation}\label{eq25}
\begin{array}{lll}
\mathrm{IMSE}_\mathrm{red,adj}^l(x_\mathrm{new})  &  =  &
\sum_{i=1}^l \sigma_{\delta^i}^2(x_\mathrm{new})
\left( 1+ \sum_{j=1}^{n_i} \frac{ \left(  \varepsilon_{\mathrm{LOO-CV},i}(x_j^i) - \hat{\rho}_{i-1}\varepsilon_{\mathrm{LOO-CV},i-1}(x_j^i) \right)^2 }{\sigma_{\mathrm{LOO-CV},i}^2 (x_j^i) - \hat{\rho}_{i-1}^2\sigma_{\mathrm{LOO-CV},i-1}^2 (x_j^i) } \right) \\
 & \times &
 \prod_{j=i}^{l-1}\rho_j^2 \prod_{m=1}^d \theta_i^m
\end{array}
\end{equation}
with the convention $\hat{\rho}_{0}=0$ and   $(\hat{\rho}_{i})_{i=1,\dots,l}$ is given by Proposition \ref{prop5}.
In equation (\ref{eq25}), the kriging variances $\sigma^2_{\delta^i}(x)$ in equation (\ref{eq23}) is replaced with the adjusted kriging variance presented in Subsection  \ref{subsection12}. We note that $\left(  \varepsilon_{\mathrm{LOO-CV},i}(x_j^i) - \hat{\rho}_{i-1}\varepsilon_{\mathrm{LOO-CV},i-1}(x_j^i) \right)^2$ is the part of the LOO-CV squared error explained by the bias $\delta^i(x)$ and $\sigma_{\mathrm{LOO-CV},i}^2 (x_j^i) - \hat{\rho}_{i-1}^2\sigma_{\mathrm{LOO-CV},i-1}^2 (x_j^i)$ is the corresponding LOO-CV predictive variance.
To adapt  the adjusted co-kriging variance in a multi-fidelity framework, we just have to replace $\mathrm{IMSE}_\mathrm{red}^l(x)$ with $\mathrm{IMSE}_\mathrm{red,adj}^l(x)$ in the algorithm \ref{algo1}.

\paragraph{$\mathbf{(q^i)_{i=1,\dots,s}}$ points at-a-time sequential co-kriging.}

In this paragraph, we propose an extension for the multi-fidelity model of the $q$ points at-a-time sequential design presented in Subsection \ref{subsection12}. Its principle is the following one. First, we select $q^l$ new points for the code $y^l(x)$ with the method presented in Subsection \ref{subsection12} ``LOO-CV based  $q$ points at-a-time Sequential design''. Then, we consider these points as known for the code $y^{l-1}(x)$ and we select $q^{l-1}$ new points for this  code with the same method. We note that, as presented in Subsection \ref{subsection11}, we can use a liar method to compute the new co-kriging variance without simulating $y^{l-1}(x)$ at the $q^l$ new points. Finally, we repeat this procedure for all code levels from  $y^{l-2}(x)$ to $y^{1}(x)$. At the end of the procedure, we have $\sum_{i=j}^{l} q^i$ new points at     level $j$ and we want to find the allocation $\{q^1,\dots,q^l\}$ leading the largest  potential uncertainty reduction and under the constraint of a constant CPU time budget.  
We note   the CPU time budget $T = \sum_{j=1}^l \sum_{i=j}^{l} q^i t^j $ where $(t^i)_{i=1,\dots,s}$ represents the CPU times of codes $(y^i(x))_{i=1,\dots,s}$. 
The  algorithm \ref{algo2} presents the suggested  $q$ points at-a-time sequential co-kriging.

\begin{algorithm}
\caption{$(q^i)_{i=1,\dots,s}$  points at-a-time sequential co-kriging}
\label{algo2}
\begin{algorithmic}[1]
\State Set the  budget $T > 0$ and the allocation $\{q^1,\dots,q^l\}$ such that $\sum_{j=1}^l \sum_{i=j}^{l} q^i t^j = T$
\State Set  $(N_{\mathrm{MCMC}}^i)_{i=1,\dots,l}$ for the M-H procedures.
\State Generate $N_{\mathrm{MCMC}}^l$ samples with respect to $k^l_{n_l}(x,x)$.
\State Find the $N^l$ cluster centers $(C_i^l)_{i=1,\dots,N^l}$ such that $N^l = \max_{N\geq q^l} \min_{x \in (C_i^l)_i}  k^l_{n_l}(x,x)$
\State Select from $(C_i^l)_{i=1,\dots,N^l}$  the $q^l$ points $(x_\mathrm{new,i}^l)_{i=1,\dots,q_l} $ where $k^l_{n_l,\mathrm{adj}}(x,x)$ is the largest.
\For { $m=l-1,\dots,1$}
\State Compute $k^m_{n_m+\sum_{i=m+1}^l q^{i}}(x,x)$ with the new  points $\left((x_\mathrm{new,i}^j)_{i=1,\dots,q_l} \right)_{j=m+1,\dots,l} $
\State Generate $N_{\mathrm{MCMC}}^m$ samples with respect to $k^m_{n_m+\sum_{i=m+1}^l q^{i}}(x,x)$.
\State Find the $N^m$ cluster centers $(C_i^m)_{i=1,\dots,N^m}$ such that $N^m = \max_{N\geq q^m } \min_{x \in (C_i^m)_i}  k^m_{n_m+\sum_{i=m+1}^l q^{i}}(x,x)$
\State Select from $(C_i^m)_{i=1,\dots,N^m}$  the $q^m$ points $(x_\mathrm{new,i}^m)_{i=1,\dots,q_m} $ where $k^m_{n_m+\sum_{i=m+1}^l q^{i}, \mathrm{adj}}(x,x)$ is the largest.
\EndFor
\end{algorithmic}
\end{algorithm}

In Algorithm \ref{algo2},  $k^l_{n_l+\sum_{i=l+1}^s q^{i}}(x,x)$ corresponds to the kernel  of the Gaussian process $Y^l_{n_l}(x)$ conditioned by the observations $\left((x_\mathrm{new,i}^j)_{i=1,\dots,q^s} \right)_{j=l+1,\dots,s} $   when the parameters $(\sigma_i^2)_{i=1,\dots,l}$ and $(\theta_i)_{i=1,\dots,l}$ are considered as known (i.e. this  corresponds to a liar method).
Furthermore,   $k^l_{n_l+\sum_{i=l+1}^s q^{i}, \mathrm{adj}}(x,x)$ corresponds to the predictor's variance $k^l_{n_l+\sum_{i=l+1}^s q^{i}}(x,x)$ adjusted with the LOO-CV errors and variances:
\begin{equation}\label{eq26}
\begin{array}{lll}
k^l_{n_l+\sum_{i=l+1}^s q^{i}, \mathrm{adj}}(x,x) &   =  &  \sum_{i=1}^l \sigma_{\delta^i+\sum_{i=l+1}^s q^{i}}^2(x )
\left( 1+ \sum_{j=1}^{n_i} \frac{  \left(  \varepsilon_{\mathrm{LOO-CV},i}(x_j^i) - \hat{\rho}_{i-1}\varepsilon_{\mathrm{LOO-CV},i-1}(x_j^i) \right)^2}{\sigma_{\mathrm{LOO-CV},i}^2 (x_j^i) - \hat{\rho}_{i-1}^2\sigma_{\mathrm{LOO-CV},i-1}^2 (x_j^i) } \right) \\
& \times &
 \prod_{j=i}^{l-1}\rho_j^2 \prod_{m=1}^d \theta_i^m
\end{array}
\end{equation}
where $k^1_{n_1+\sum_{i=l+1}^s q^{i}}$ and $\sigma_{\delta^i+\sum_{i=l+1}^s q^{i}}^2(x_\mathrm{new})$ are deduced from the equation  (\ref{eq13}).
We note that for the M-H procedures, we use a Gaussian jumping distribution with a standard deviation such that acceptance rate is around 30\%.

Furthermore, let us consider the following quantity
\begin{equation}\label{eq27}
\mathrm{IMSE}_\mathrm{red,q} = \sum_{i=1}^l \sum_{r=1,\dots,q^i} \sigma_{\delta^i}^2(x_{\mathrm{new},r}^i) \prod_{j=i}^{l-1}\rho_j^2 \prod_{m=1}^d \theta_i^m
\end{equation}
We consider the allocation $\{q^1,\dots,q^l\}$ which solves the following optimization problem:
\begin{equation}\label{eq28}
\{q^1,\dots,q^l\} = \arg \max_{\{q^1,\dots,q^l\}}  \mathrm{IMSE}_\mathrm{red,q} \, \, \mathrm{such \, that} \, \, \sum_{j=1}^l \sum_{i=j}^{l} q^i t^j = T
\end{equation}
i.e. we look for the allocation leading the maximal  uncertainty reduction.
This optimization problem is very complex to solve and a sub-optimal solution will be often considered. Nevertheless, when the number of code levels and the budget $T$ are  low (e.g. $s=2$ in our application)  an exhaustive exploration of the allocation $\{q^1,\dots,q^l\}$ can be performed. We are in that case in the presented application . Furthermore, we note that $\mathrm{IMSE}_\mathrm{red,q}$ is a proxy on the IMSE reduction when we add $\left( (x_\mathrm{new,i}^m)_{i=1,\dots,q_m} \right)_{m=1,\dots,l} $ at code levels $(y^m(x))_{m=1,\dots,l}$.

In practical application, the algorithm \ref{algo2} is reiterated until we reach a prescribed precision or the computational time budget is exhausted. 

\section{Applications}

We compare in this Section the MinIMSE, KleiCrit and AdjMMSE criteria on toy examples and on an application concerning a spherical tank under pressure. We present both the cases of 1 point at-a-time and $q$ points at-a-time sequential kriging. Then, we compare on the tank application, the suggested sequential kriging and co-kriging methods with $s=2$ levels.  The purpose of this section is to emphasize the efficiency of the LOO-CV-based criteria and to highlight  that a multi-fidelity  analysis can be worthwhile. Finally, for the multi-fidelity sequential co-kriging, we present the allocation of the simulations between the coarse code and the accurate one. We note that for the different examples, we compare the different methods given a prescribed computational time budget.

\subsection{Comparison between sequential kriging criteria}

In this subsection, the 1 point at-a-time sequential kriging criteria (MinIMSE, KleiCrit, AdjMMSE) are compared on three tabulated functions:
\begin{itemize}
\item Ackley's function on $[-2,2]^2$ \cite{Ack87}:
\end{itemize}
\begin{displaymath}
f(x,y) = -20 \mathrm{exp}\left( -0.2 \sqrt{\frac{x^2+y^2}{2}} \right) - \mathrm{exp}\left( \frac{\mathrm{cos}(2\pi x)  + \mathrm{cos}(2\pi y)}{2} \right) + 20 + \mathrm{exp}(1)
\end{displaymath}
\begin{itemize}
\item Shubert's function on $[-2,2]^2$ \cite{XianLi01}:
\end{itemize}
\begin{displaymath}
f(x,y) = \left(\sum_{k=1}^5 {k \mathrm{cos}\left((k+1)x + k  \right)} \right) \left( \sum_{k=1}^5 {k \mathrm{cos}\left((k+1)y  + k  \right)} \right)
\end{displaymath}
\begin{itemize}
\item Michalewicz's function on $[0,\pi]^2$ \cite{Mich92}:
\end{itemize}
\begin{displaymath}
f(x,y) = - \mathrm{sin}\left( x\right) \left( \mathrm{\sin}\left( \frac{x^2}{\pi} \right) \right)^{20}- \mathrm{sin}\left( y\right) \left( \mathrm{\sin}\left( \frac{y^2}{\pi} \right) \right)^{20}
\end{displaymath}
The comparison is performed on a test set $\mathbf{D}_\mathrm{test}$ composed of $n_\mathrm{test} = 1000$ points uniformly spread on the input parameter space and from 50 different initial experimental design sets. We compare the different methods with respect to the Normalized RMSE:
\begin{equation}
\mathrm{Norm \, RMSE}  = \frac{\sqrt{\sum_{i=1}^{n_\mathrm{test}} \left( y_\mathrm{real}(x^i_\mathrm{test}) - y_\mathrm{pred}(x) \right)^2/n_\mathrm{test}}}{\max_{x \in \mathbf{D}_\mathrm{test}} y_\mathrm{real}(x) - \min_{x \in \mathbf{D}_\mathrm{test}} y_\mathrm{real}(x) }
\end{equation}
where $y_\mathrm{real}(x) $ is the real value of the output and $y_\mathrm{pred}(x) $ the predicted one.
The initial experimental design sets  are LHS designs of 10 points   optimized with respect to the S-optimality \cite{Sto05}. From these designs,  50 sequential krigings are performed and the convergence of the mean and the quantiles of the Normalized RMSE are computed for the three criteria. Furthermore, after each added point, the parameters of the kriging model are re-estimated with a maximum likelihood method. These estimations are performed thanks to the R library 'DiceKriging' \cite{DiceKriging}. 

\begin{figure}[ht]
\begin{center}
\vskip-04ex
\includegraphics[trim = 0cm 0.5cm 0cm 0cm , width = 0.45\linewidth]{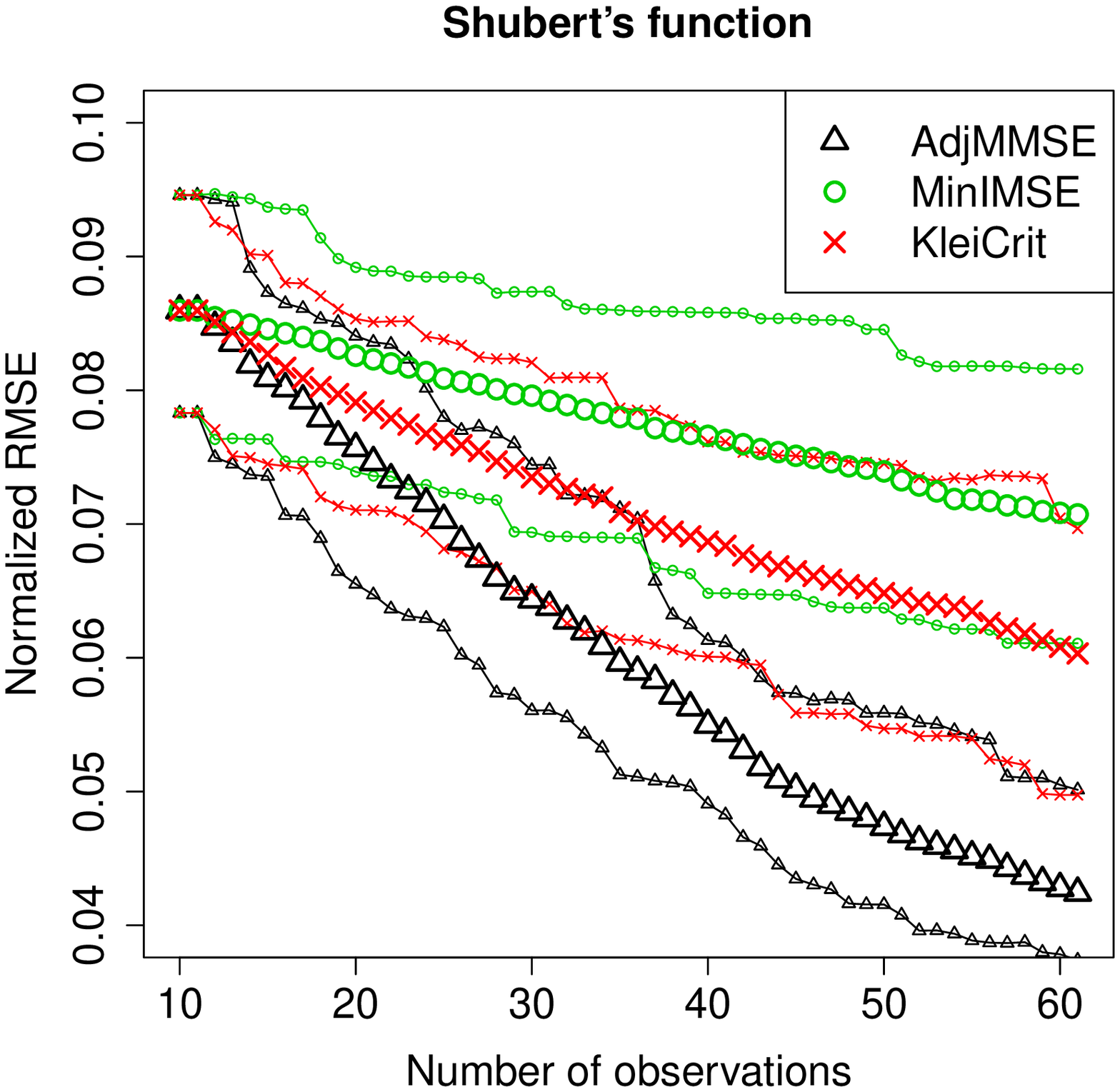}
\includegraphics[trim = 0cm 0.5cm 0cm 0cm , width = 0.45\linewidth]{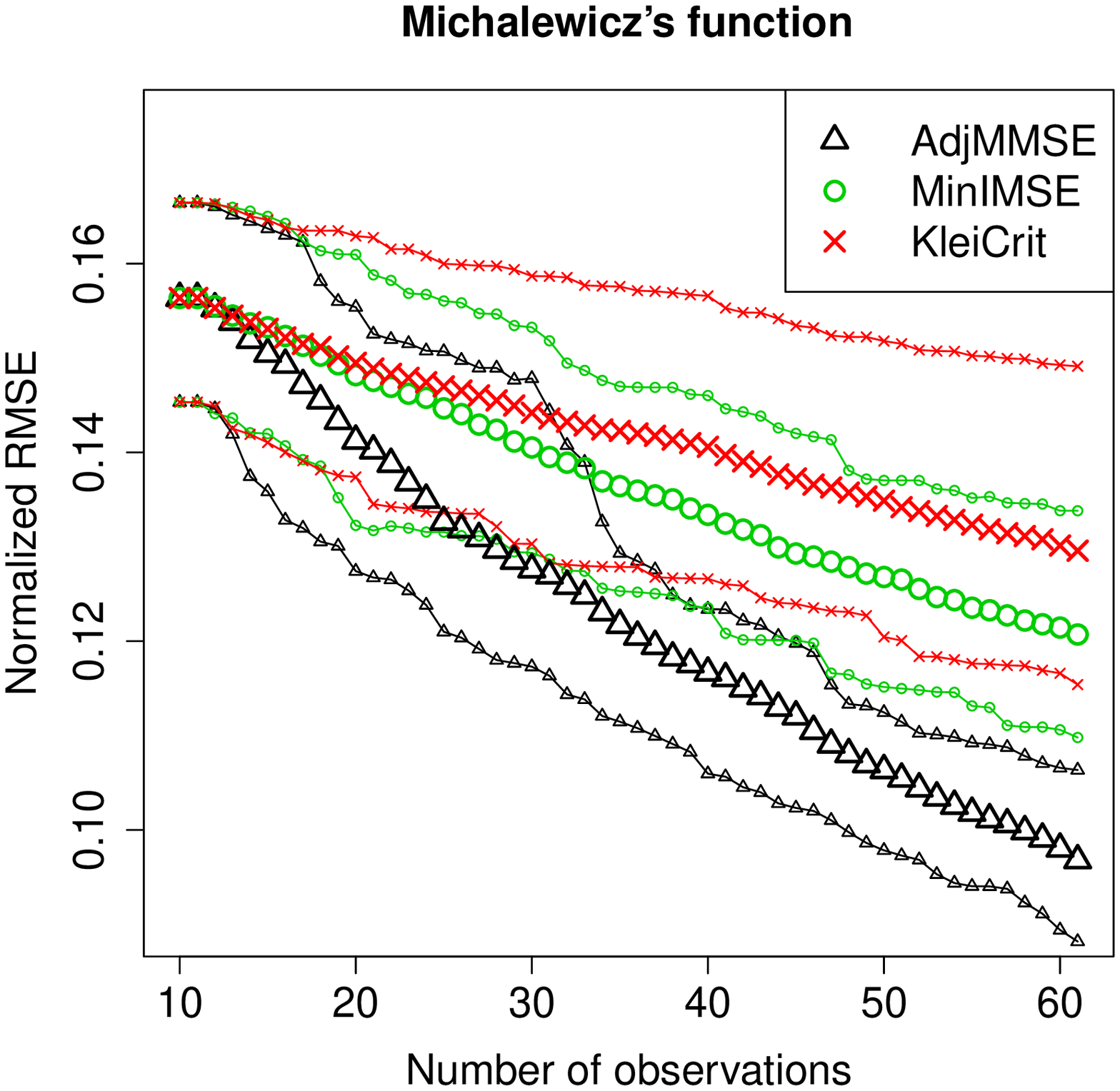}
\includegraphics[trim = 0cm 0.5cm  0cm 0cm , width = 0.45\linewidth]{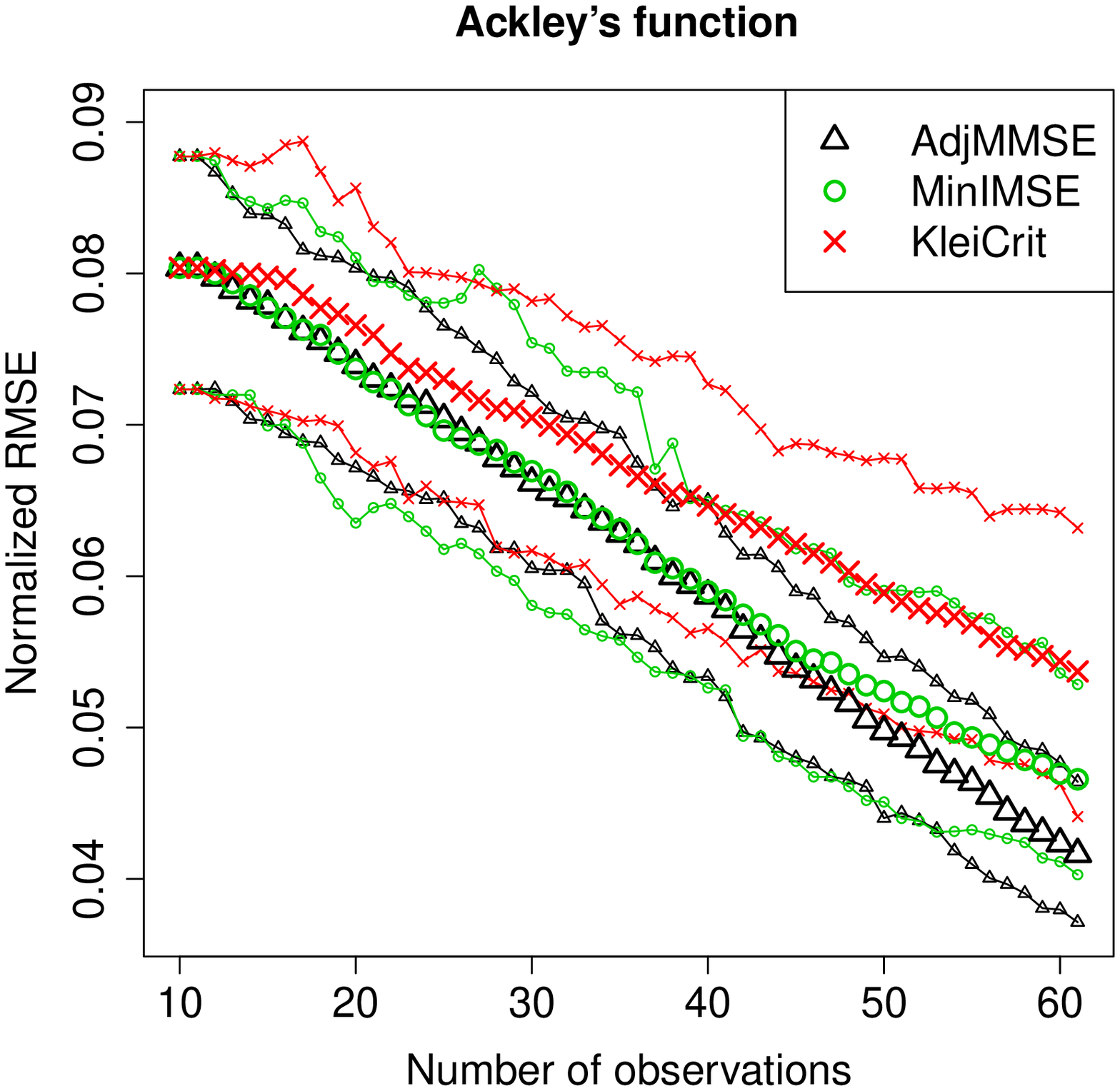}
\vskip-03ex
\caption{Comparison between  1 point at-a-time sequential kriging criteria on toy examples. The bold triangles represent the mean of the empirical MSE for the AdjMMSE criterion, the bold circles represent it for the MinIMSE criterion and the bold Crosses represent it for the KleiCrit criterion. Furthermore, the thin triangles, circles and crosses represent the quantiles of probabilities 10\% and 90\% of the empirical MSE.}
\label{Toy_compkm}
\end{center}
\end{figure}

Figure \ref{Toy_compkm} illustrates the efficiency of the   criterion AdjMMSE. Indeed, for the Shubert's and Michalewicz's functions, we see that the accuracy of the 1 point at-a-time kriging with this criterion is significantly better than the one of the others criteria (both in terms of mean and quantiles of the  Normalized RMSE). In fact, these functions have the particularity to have non-stationarity   in some areas of the input parameter space. Thus, the errors are more  important in these locations  and the suggested criterion focus the new points  on it. Furthermore, the non-stationarity is particularly  important for  the Shubert's function. For this reason, the IMSE criterion performed very bad in that case. Indeed, this criterion is   efficient for stationary functions (i.e. when the predictor's MSE well predicts the real model errors). In contrast,  the jackknife predictor's MSE provided by the criterion KleiCrit manages to catch this non-stationarity and it performs better than the IMSE criterion.
Moreover, we see  that the performance of the AdjMMSE and IMSE criteria are equivalent for  the Ackley's function. We note that the Ackley's function  is perfectly stationary and the errors have thus the same order of magnitude over the input parameter space. 

These examples  illustrate  the fact that our criterion is more efficient than the other criteria when the functions are non-stationary and it  remains efficient even in the cases where the functions are stationary (its efficiency is equivalent to  the one of the IMSE criterion).

\subsection{Spherical tank under internal pressure example}

\begin{figure}[ht]
\begin{center}
\vskip-04ex
\includegraphics[trim = 7cm 3cm 1cm 1cm , width = 0.30\linewidth]{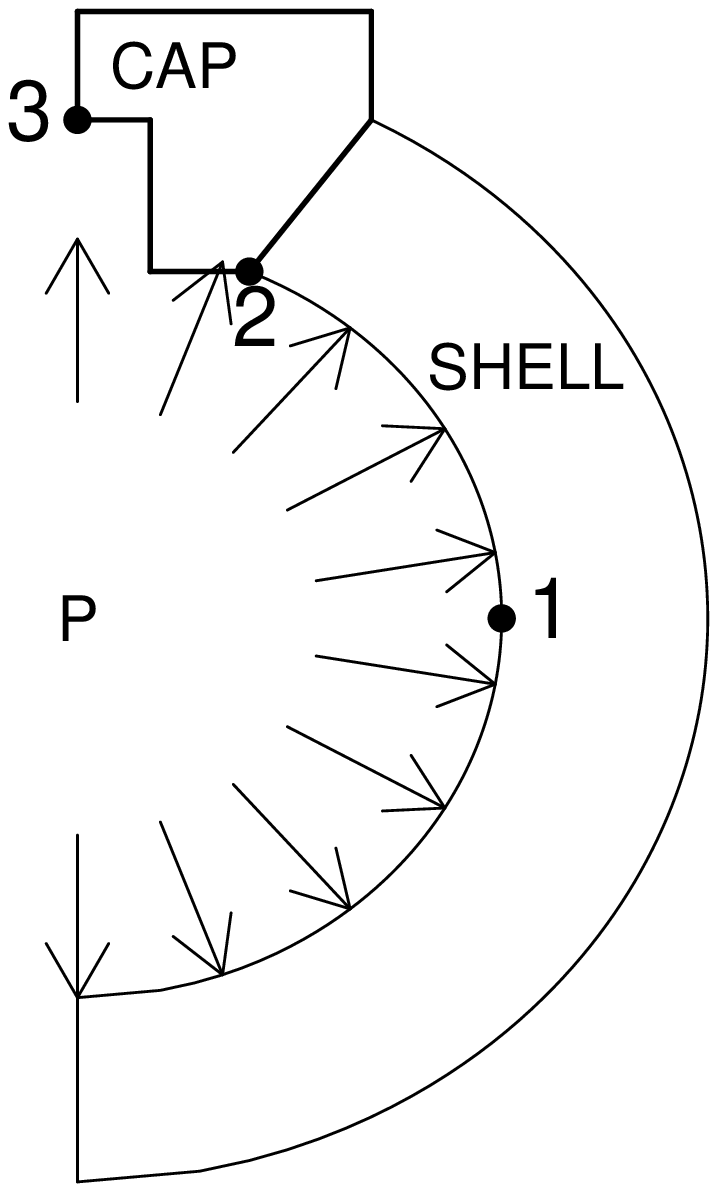}
\vskip-03ex
\caption{Scheme of the  spherical tank under pressure. }
\label{scheme}
\end{center}
\end{figure}

In this section, we deal  with an example about a spherical tank under internal pressure.  We are interested in the von Mises stresses on the  three   points labeled  in figure  \ref{scheme}. 
Indeed, we want to prevent from material yielding which occurs when the von Mises stress reaches the critical  yield strength.

The system illustrated in figure \ref{scheme} depends on the following parameters:
\begin{itemize}
\item $P \, (MPa) \in [30, 50]$: the value of the  internal pressure.
\item $R_{int} \, (mm) \in [1500, 2500]$:  the length of the internal  radius of the shell.
\item $T_{shell} \, (mm) \in [300, 500]$: the thickness of the shell.
\item $T_{cap} \, (mm) \in [100, 300]$: the thickness of the cap.
\item $E_{shell} \, (GPa) \in [63, 77]$: the Young's modulus of the shell material.
\item $E_{cap} \, (GPa) \in [189, 231]$: the Young's modulus of the cap material.
\item $\sigma_{y,shell} \, (MPa) \in [200, 300]$: the yield stress of the cap material.
\item $\sigma_{y,cap} \, (MPa) \in [400, 800]$: the yield stress of the cap material.
\end{itemize}
The accurate code output $y^2(x)$ is the value of the von Mises stress provided by an Aster  finite elements code (\url{http://www.code-aster.org}) modelling the system presented in figure \ref{scheme}. We use the notation $x = (P, R_{int}, T_{shell}, T_{cap}, E_{shell}, E_{cap}, \sigma_{y,shell}, \sigma_{y,cap})$. We note  that the material properties of the shell correspond to high quality aluminum and the ones  of the cap corresponds to steel from classical   to  high quality.  Then, the coarse code output $y^1(x)$  is the value of the von Mises stress given by the 1D simplification of the tank  (\ref{eq27y})  (it corresponds to a perfect spherical tank under pressure, i.e. without cap).
\begin{equation}\label{eq27y}
y^1(x) = \frac{3}{2} \frac{\left( R_{int} + T_{shell} \right)^3}{\left( R_{int} + T_{shell} \right)^3 - R_{int}^3}P
\end{equation}
According to equation (\ref{eq27y}), the actual input dimension of $y^1(x)$ is three (it depends only on $P$, $R_{int}$ and $T_{shell}$) while  a sensitivity analysis performed with a Sobol decomposition gives   that the accurate code depends essentially on four parameters ($P$, $R_{int}$, $T_{shell}$ and $T_{cap}$). Furthermore, the response is highly stationary. Therefore, only few points are necessary to well predict the output of the code. For these reasons, we can start the sequential strategies from an initial experimental design set with only 10 points. 

Thus, for the different comparisons, we use a  S-optimal  LHS design $\mathbf{D}^2$ of 10 points  for the code $y^2(x)$. For the coarse code $y^1(x)$, we start with a design $\mathbf{D}^1$  of 20 points.
It is created with the following  procedure. First, we create a S-Optimal design $\mathbf{\tilde{D}}^1$ of 20 points. Second, we remove from $\mathbf{\tilde{D}}^1$ the 10 points that are the closest to those of $\mathbf{D}^2$. Finally, $\mathbf{D}^1$ is the concatenation of $\mathbf{D}^2$ and $\mathbf{\tilde{D}}^1$ (this procedure ensures the nested property $\mathbf{D}^2 \subset \mathbf{D}^1 $). We note that   the CPU time is around 1 minute  for the accurate code and $10^{-8}$ seconds for the coarse code. Nevertheless, to be in a more realistic case, we consider that the CPU time ratio  between $y^2(x)$ and $y^1(x)$ equals $B_{2/1} = 10$. Furthermore, each sequential procedure is performed with 40 different initial design sets. Then, the mean and the quantiles of probabilities 90\% and 10\% of the empirical Normalized MSE are computed from a test set composed of 1000 points uniformly spread on the input parameter space. Finally, for the M-H  procedure, we use a Gaussian jumping distribution such that the acceptance rate is around 30\% and we set   $N_{\mathrm{MCMC}} = 50 000$ (we use 5 000 samples for the the burn-in procedure of the M-H method, see \cite{RobCas04}). For the M-H  procedure, we use the package R CRAN mcmc. We note that after each added points, the parameters of the kriging or co-kriging models are re-estimated with a maximum likelihood method.

The remainder of this section is organized as follows. First we compare the MSE of the 1 point at-a-time sequential kriging with the one of the $q=5$ points at-a-time one. Second, we compare for a given CPU time budget the sequential kriging and cokriging strategies. In the forthcoming developments, the response $i=1, 2, 3$  refers to the value of the von Mises stress at point $i$ on figure \ref{scheme}.

\subsubsection{Comparison between sequential kriging criteria}

\begin{figure}[ht]
\begin{center}
\vskip-04ex
\includegraphics[trim = 0cm 0.5cm 0cm 0cm , width = 0.45\linewidth]{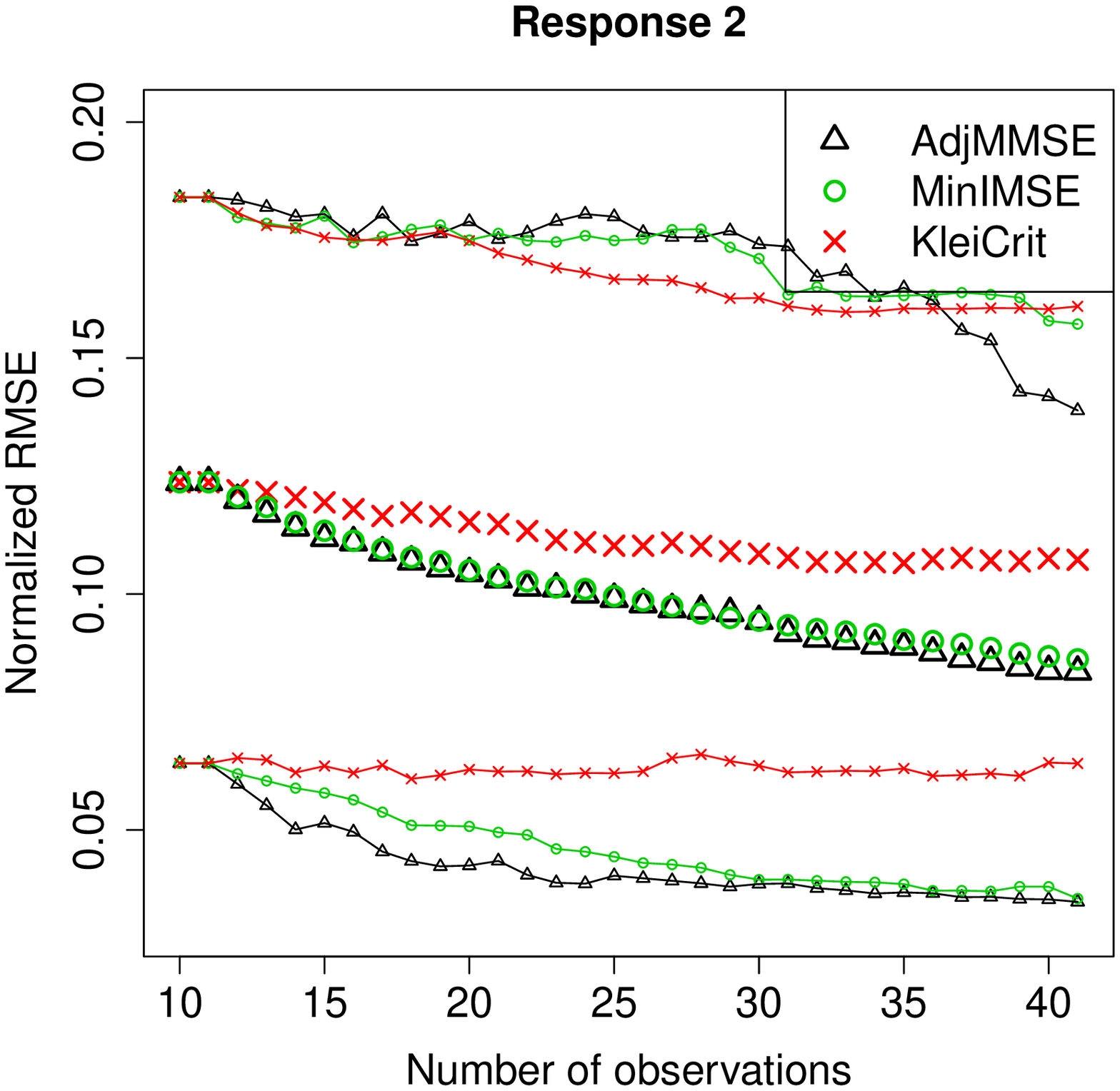}
\includegraphics[trim = 0cm 0.5cm 0cm 0cm , width = 0.45\linewidth]{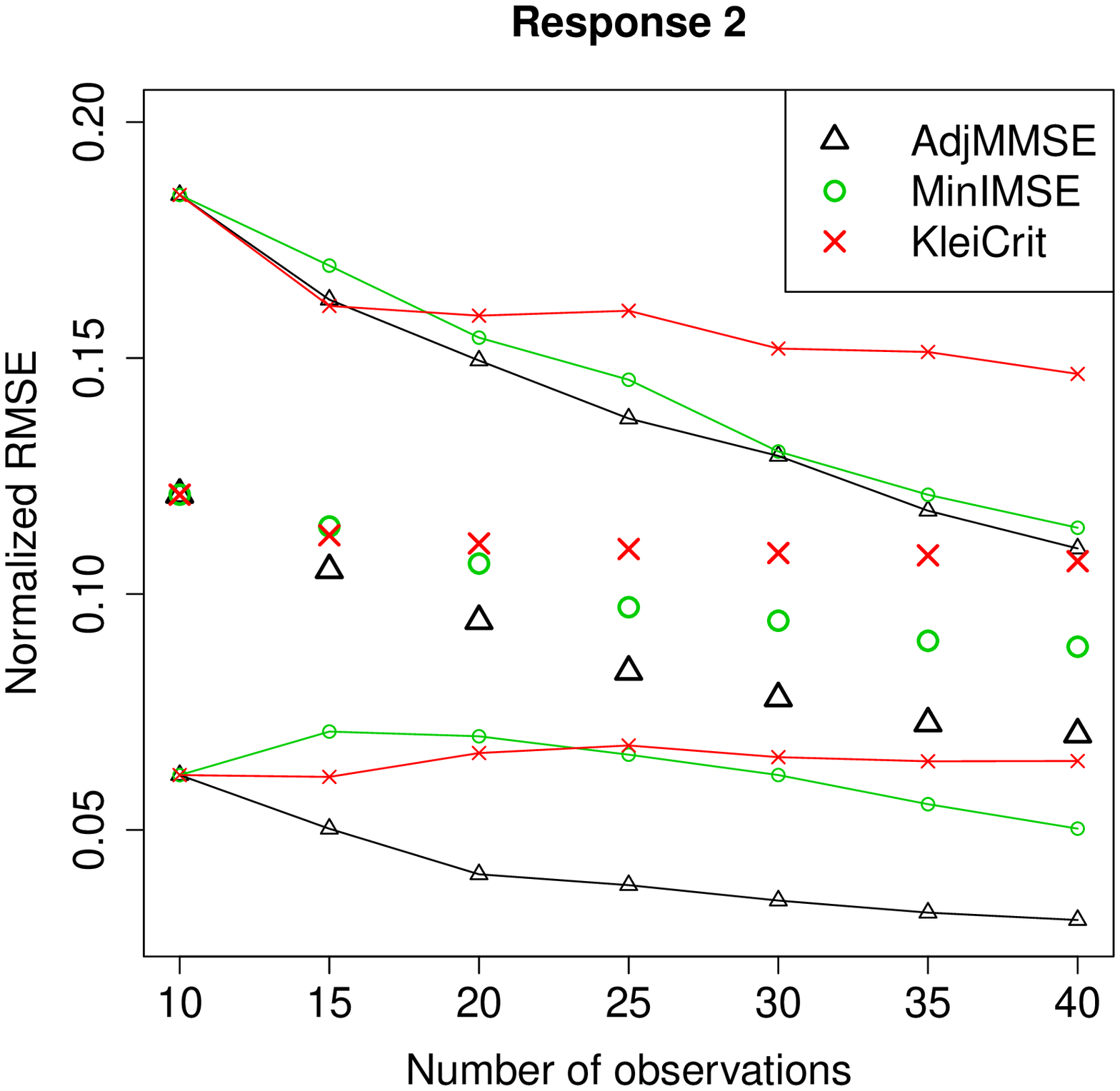}
\vskip-03ex
\caption{Comparison between  1 point at-a-time sequential kriging criteria (on  left) and $q = 5$ points  at-a-time sequential kriging criteria (on right) on the spherical tank example. The bold triangles represent the mean of the empirical MSE for the AdjMMSE criterion, the bold circles represent it for the MinIMSE criterion and the bold Crosses represent it for the KleiCrit criterion. Furthermore, the thin triangles, circles and crosses represent the quantiles of probabilities 10\% and 90\% of the empirical MSE. }
\label{Rep3_compkm}
\end{center}
\end{figure}

Figure \ref{Rep3_compkm} compares the different criteria of  the 1 point at-a-time and the $q = 5$ points at-a-time sequential kriging. We see that the criteria MinIMSE and AdjMMSE give equivalent values for the MSE  for the 1 point at-a-time procedure and they  perform  better than the KleiCrit criterion. There are   equivalent  since the output $y^2(x)$ is perfectly stationary. Nevertheless, the criterion AdjMMSE is the most efficient for the $q=5$ points at-a-time procedure. Indeed, the 5 points provided by a liar method with the MinIMSE criterion are not necessarily  those  which maximize the reduction of the IMSE. The  method suggested in section \ref{subsection23}  seems to give  better solution.

\subsubsection{Comparison between kriging and co-kriging sequential analysis}

\begin{figure}[ht]
\begin{center}
\vskip-04ex
\includegraphics[trim = 0cm 0.5cm 0cm 0cm , width = 0.45\linewidth]{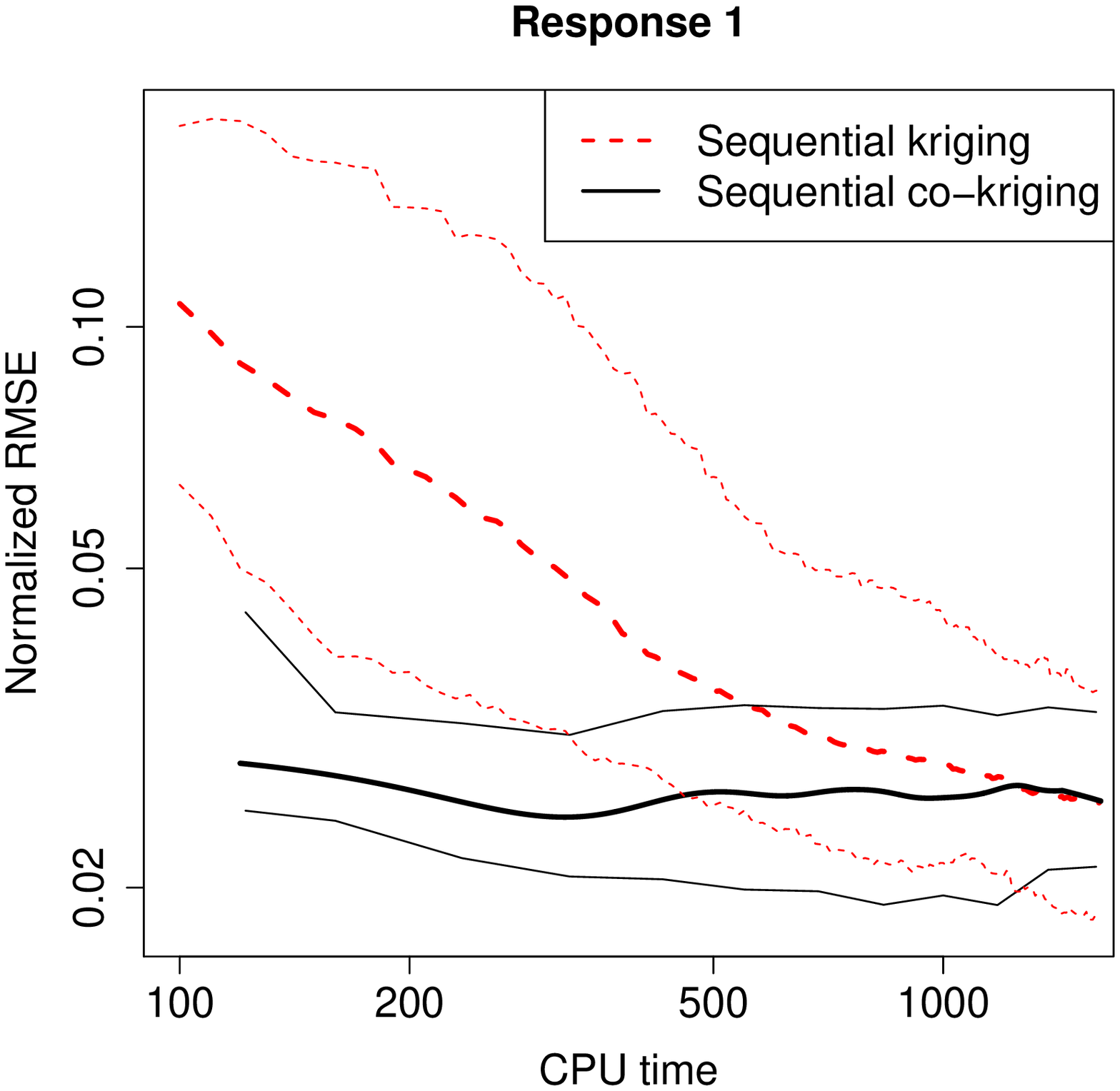}
\includegraphics[trim = 0cm 0.5cm 0cm 0cm , width = 0.45\linewidth]{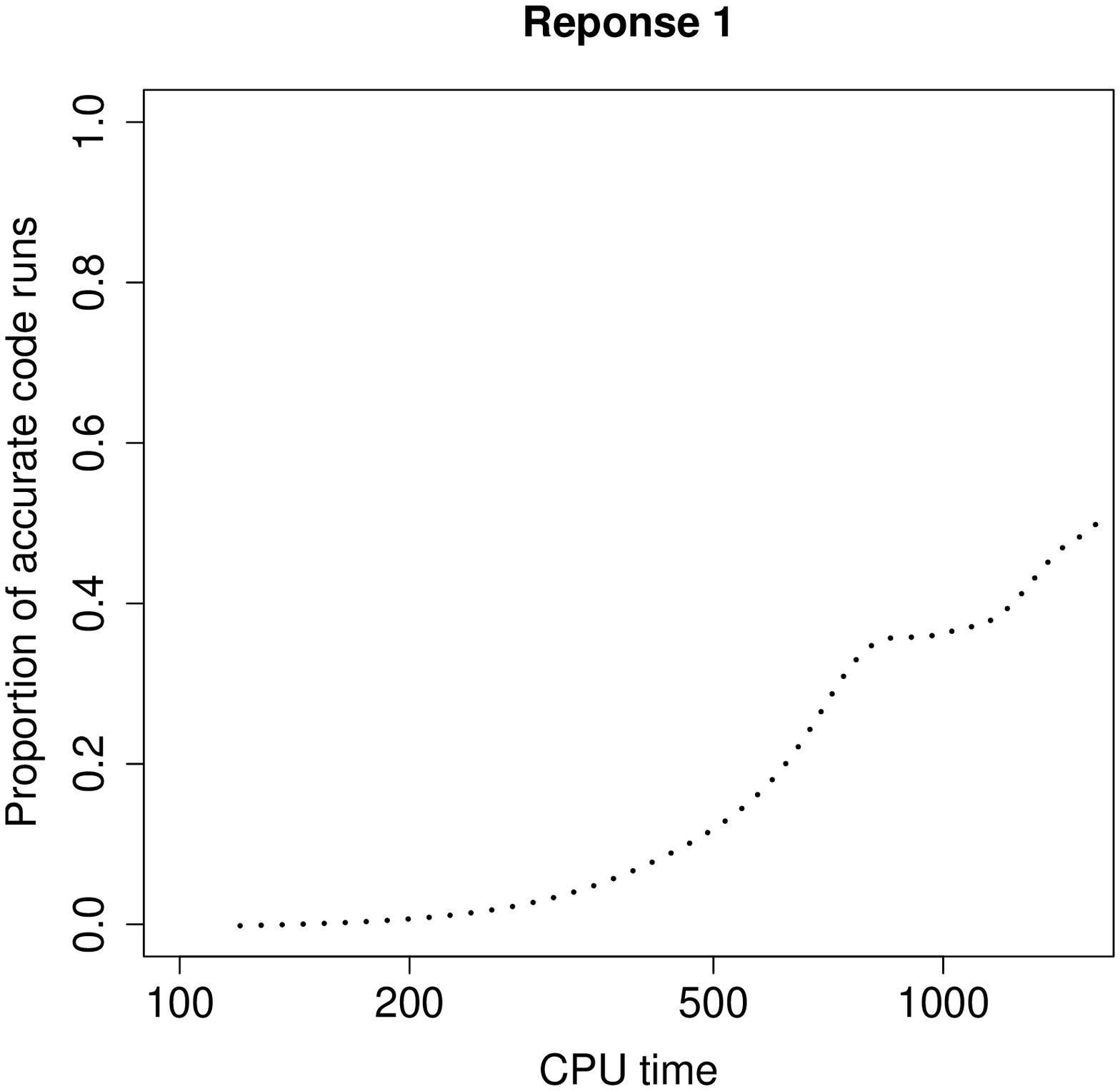}
\vskip-03ex
\caption{Comparison between  1 point at-a-time   sequential kriging and co-kriging   on the response 1 of the spherical tank example with respect to the AdjMMSE criterion (on the left). The thick dashed line represents  the mean of the empirical MSE for the   sequential kriging and the thick solid line represents  it for the sequential co-kriging. The thin   lines represent the quantiles of probabilities 10\% and 90\% of the empirical MSE. Figure on the right represents the proportion of runs allocated to the accurate code.  }
\label{Rep1_kmvscokm}
\end{center}
\end{figure}

In this section, we compare the  sequential kriging strategy with the sequential co-kriging  with  respect to the AdjMMSE criterion.
Figure \ref{Rep1_kmvscokm} gives the convergence of the empirical normalized MSE for the response 1. We see that the sequential co-kriging performs  better than the kriging one. Furthermore, at the beginning of the method,   the proportion of runs for the accurate  code is very low. Indeed, the coarse code and the accurate code are extremely correlated for this response (around 99\%) and thus, during the sequential strategy,  the bias between the two codes is well estimated. Then, when the coarse code is well approximated, the sequential strategy starts to  run the accurate one  (for a CPU time around 500).

\begin{figure}[ht]
\begin{center}
\vskip-04ex
\includegraphics[trim = 0cm 0.5cm 0cm 0cm , width = 0.45\linewidth]{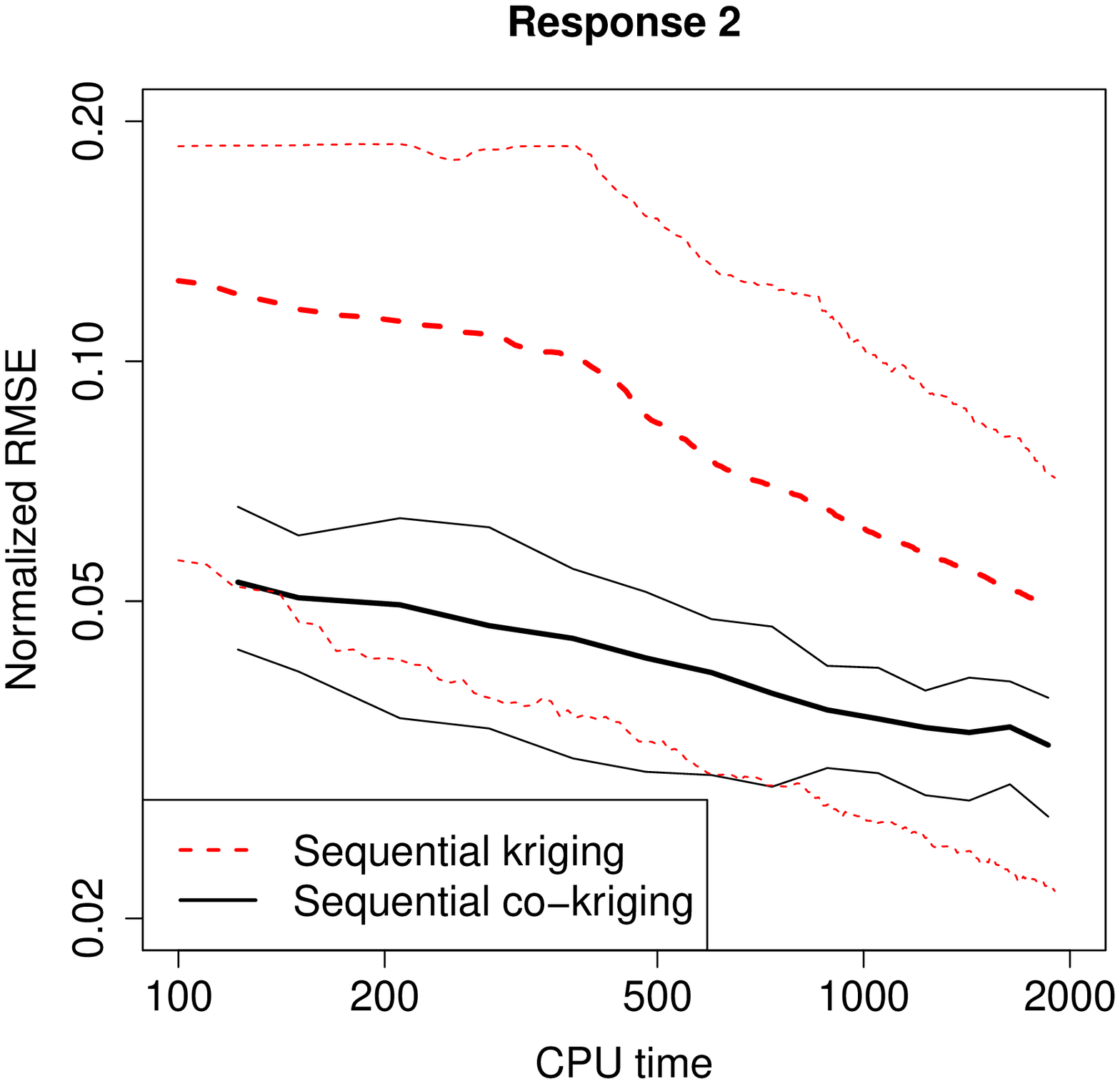}
\includegraphics[trim = 0cm 0.5cm 0cm 0cm , width = 0.45\linewidth]{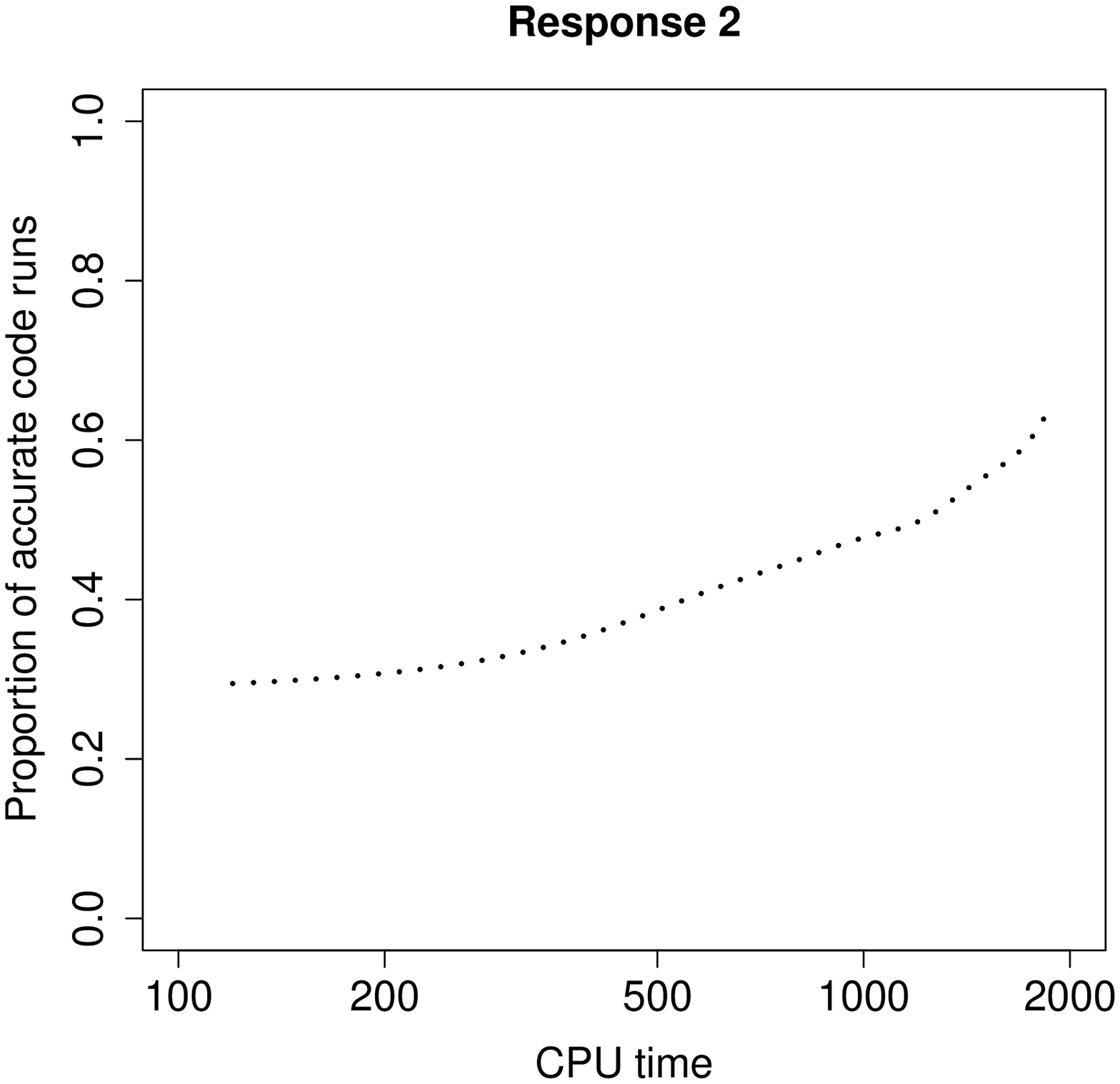}
\vskip-03ex
\caption{ Comparison between  1 point at-a-time   sequential kriging and co-kriging   on the response 2 of the spherical tank example with respect to the AdjMMSE criterion (on the left). The thick dashed line represents  the mean of the empirical MSE for the   sequential kriging and the thick solid line represents  it for the sequential co-kriging. The thin   lines represent the quantiles of probabilities 10\% and 90\% of the empirical MSE.  Figure on the right represents the proportion of runs allocated to the accurate code.}
\label{Rep3_kmvscokm}
\end{center}
\end{figure}

Figure \ref{Rep3_kmvscokm} gives the convergence of the errors for the response 2. For this response, the correlation between the coarse and the accurate code is around 80\%. Therefore, the proportion of runs for the accurate code determined by the sequential strategy  is more important than in Figure \ref{Rep1_kmvscokm}. Furthermore, we see that this proportion increases with the CPU time. It means that the sequential co-kriging improves the approximation of the coarse code at the beginning of the procedure and then focuses on the accurate code. As a result, we see that the sequential co-kriging strategy is substantially better than the kriging one.

\begin{figure}[ht]
\begin{center}
\vskip-04ex
\includegraphics[trim = 0cm 0.5cm 0cm 0cm , width = 0.45\linewidth]{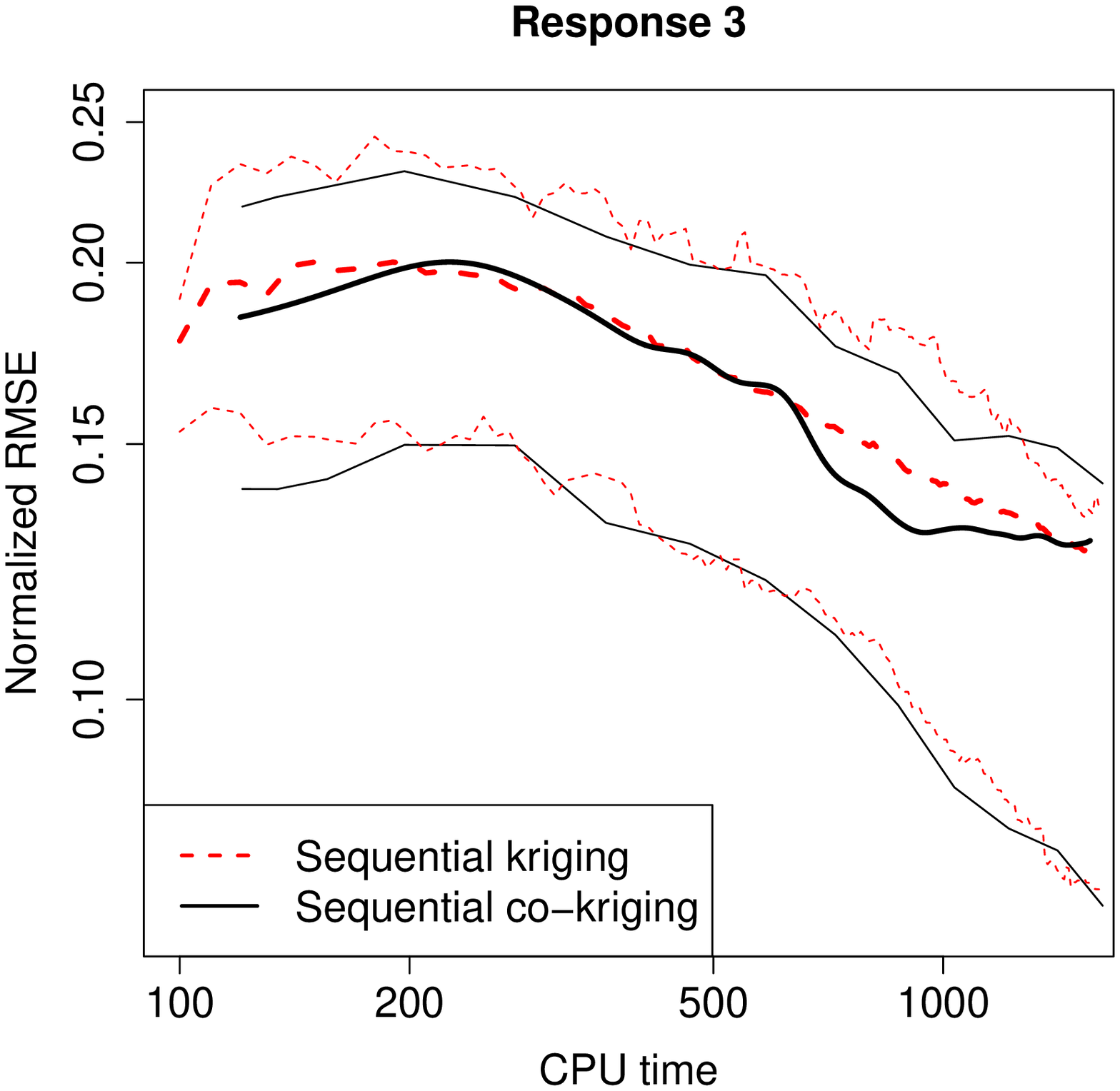}
\includegraphics[trim = 0cm 0.5cm 0cm 0cm , width = 0.45\linewidth]{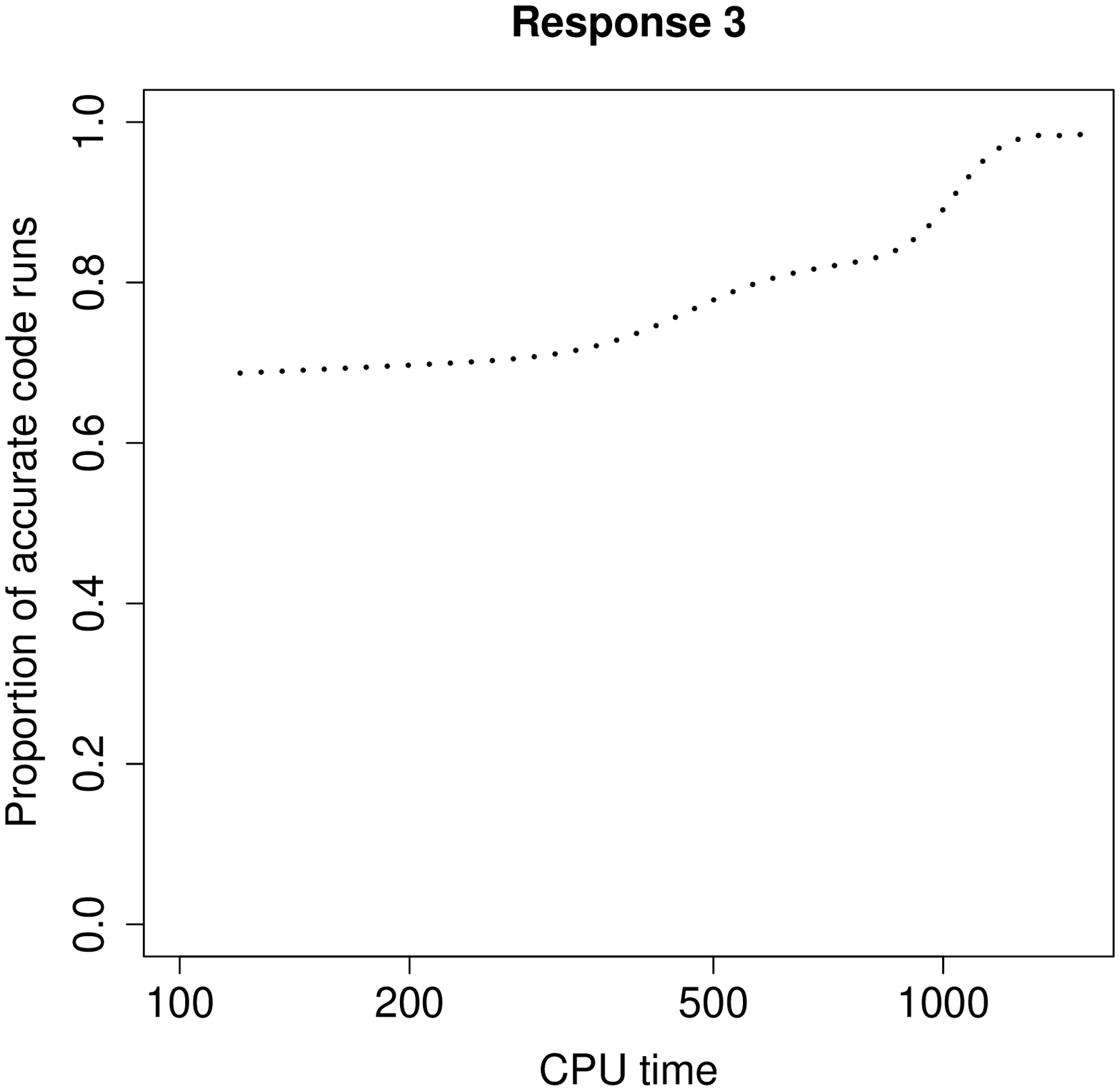}
\vskip-03ex
\caption{Comparison between  1 point at-a-time   sequential kriging and co-kriging   on the response 3 of the spherical tank example with respect to the AdjMMSE criterion (on the left). The thick dashed line represents  the mean of the empirical MSE for the   sequential kriging and the thick solid line represents  it for the sequential co-kriging. The thin   lines represent the quantiles of probabilities 10\% and 90\% of the empirical MSE.  Figure on the right represents the proportion of runs allocated to the accurate code.  }
\label{Rep4_kmvscokm}
\end{center}
\end{figure}

Figures  \ref{Rep1_kmvscokm} and  \ref{Rep3_kmvscokm} illustrate the efficiency of the sequential co-kriging when the coarse code bring information on the accurate code. For the response 3, the coarse code is weakly correlated with the accurate code (around 45\%). This is due to the fact that the coarse code models the von Mises stress in  a perfect spherical tank whereas the response 3 corresponds to the  one in the cap. Figure \ref{Rep4_kmvscokm} shows that in this case, the sequential co-kriging model manages to determine  that the coarse code  is not worth being simulated. Indeed, the proportion of runs for the accurate  code is very low. Furthermore, it shows that the co-kriging sequential design performs as well as the kriging one when the coarse code is non-informative.

\begin{figure}[ht]
\begin{center}
\vskip-04ex
\includegraphics[trim = 0cm 0.5cm 0cm 0cm , width = 0.45\linewidth]{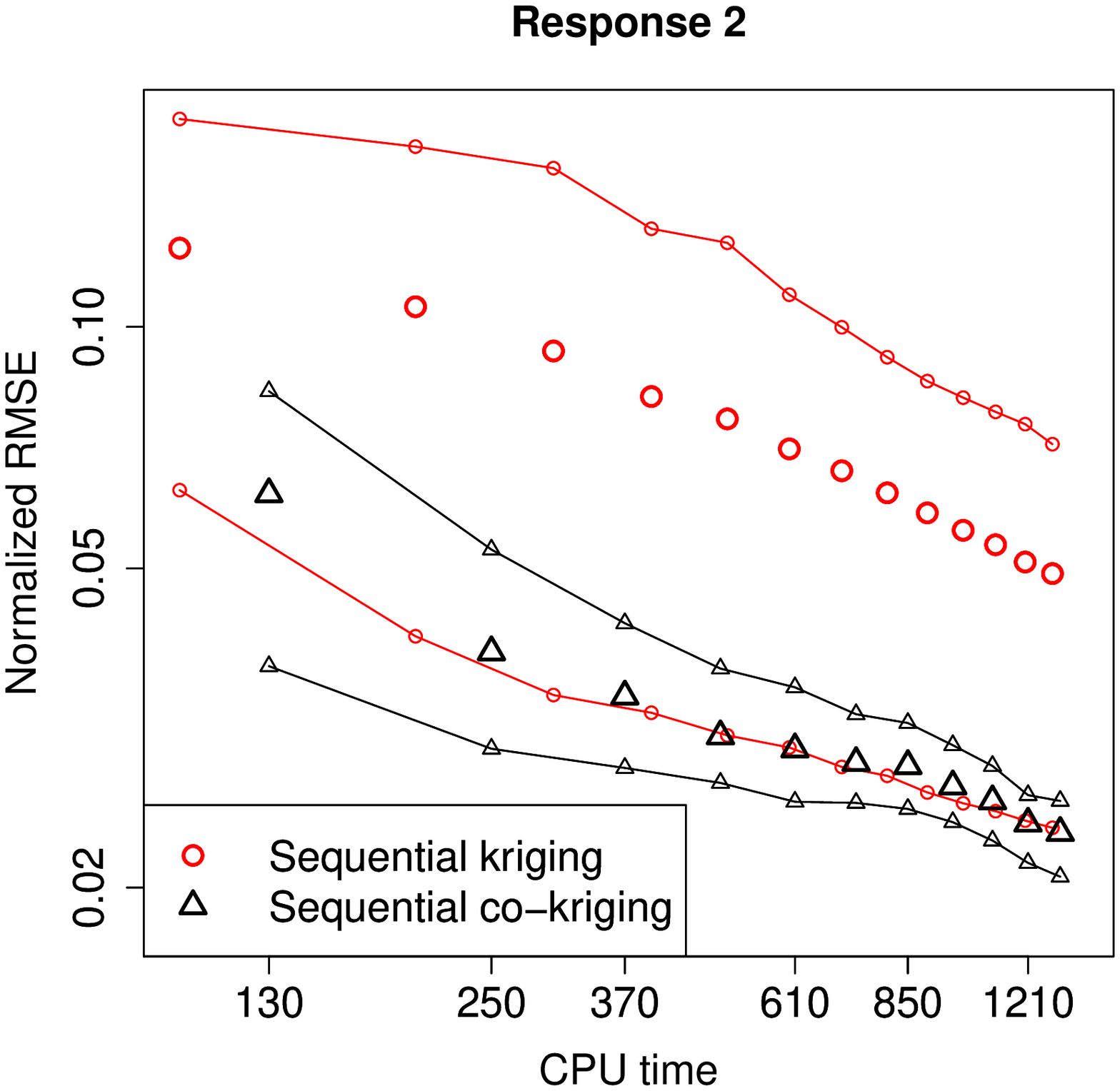}
\includegraphics[trim = 0cm 0.5cm 0cm 0cm ,width = 0.45\linewidth]{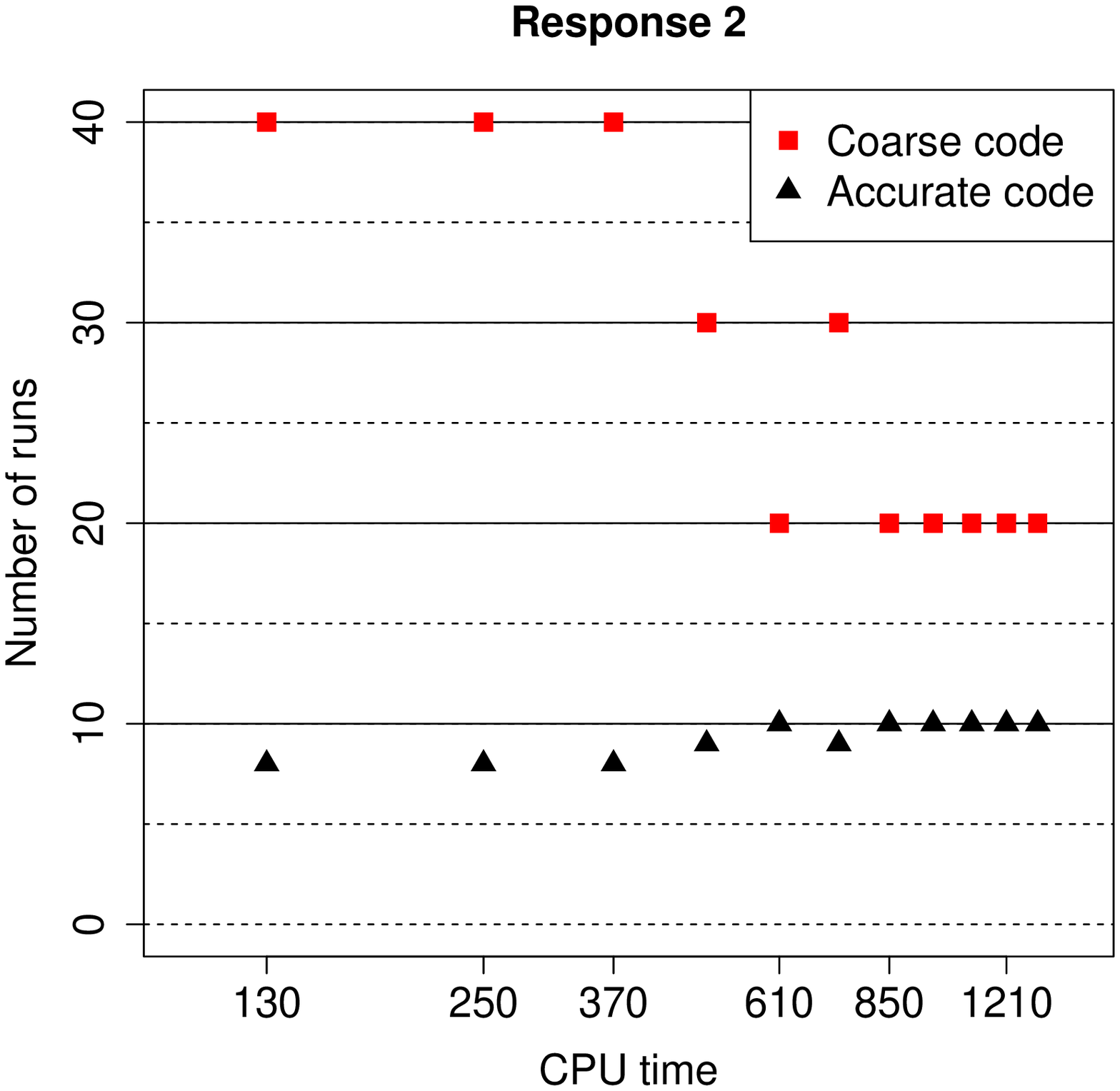}
\vskip-03ex
\caption{Comparison between $q = 10$ points at-a-time sequential kriging and $(q^1,q^2)$ points at-a-time sequential co-kriging. On the left, the bold circles represents the mean of the empirical MSE for the sequential kriging and  the bold triangles represent the one of the sequential co-kriging. Furthermore, the thin circles and triangles represent the quantiles of probabilities 10\% and 90\% of the empirical MSE. On the right, the squares represent the median number of runs for the coarse code   during  the sequential co-kriging  and   the triangles represent it for the accurate code.}
\label{Rep3_QAAT_kmvscokm}
\end{center}
\end{figure}

Finally, Figure \ref{Rep3_QAAT_kmvscokm} shows the efficiency of the $(q^1,q^2)$ at-a-time sequential co-kriging. We set in the algorithm \ref{algo2} that  $T = q^1 + q^2 + 10 q^2 = 120$ where the CPU time of the coarse code is $1$ and the one of the accurate code is $10$. For the the sequential kriging, we use a $q=10$ at-a-time sequential procedure. Furthermore,  Figure \ref{Rep3_QAAT_kmvscokm} shows  that at the beginning of the procedure, the sequential co-kriging focuses on  the approximation of the coarse code whereas at the end the number of runs for the accurate code is maximal. We note that the allocation of runs for the accurate code in figure \ref{Rep3_QAAT_kmvscokm} agrees with the proportion of runs given in Figure \ref{Rep3_kmvscokm}.

The results  of the sequential co-kriging on the different responses show that the criterion suggested in section \ref{subsection21} performs very well. Indeed, it is always better than the sequential kriging when the coarse code is informative and its performance is equivalent to it when the coarse code is not useful. Furthermore, the different proportions of runs for the accurate code emphasizes  that the criterion accurately determines  the contribution of each code to the total model error and the optimal run allocation between the  accurate and the coarse codes.

\section{Conclusion}

This paper deals with   sequential strategies for kriging and co-kriging models. The kriging model is used to approximate the output of a complex computer code and the co-kriging one allows   to improve this approximation thanks to  coarse versions of the code. 

First, we have  presented   classical sequential criteria for the kriging model and we have  suggested  another criterion based on the Leave-One-Out cross validation errors. This criterion has  allowed us to set  the new observations at locations where the model error  is  important. The examples presented in the last section have highlighted  the efficiency of the suggested criterion. Indeed, for non-stationary functions, it provides results significantly better than classical criteria  and for stationary ones its performance is equivalent to them. We have also emphasized  the performance of the suggested criterion on a real application. Furthermore, we show  in the application that when the simulations can be performed in parallel, our method has  performed better.

Second, we have presented the extension of our criterion to  multi-fidelity co-kriging models.  We have   shown in the application  that performing   a multi-fidelity sequential co-kriging is worthwhile when the coarse code versions are informative (i.e. highly correlated with the accurate code). Furthermore, a strength of the proposed approach is that it performs as well as a sequential kriging when the coarse code versions are not informative. In fact, the proposed extension takes into account the contribution of each code to the total predictor's mean squared errors and  it determines  the best run  allocation between accurate and coarse code versions  given a CPU time budget.

\section{Acknowledgements}

The author particularly thanks Professor Josselin Garnier for supervising his  work and for his fruitful guidance. He is also grateful to Dr. Gilles Defaux   for providing the  application.

\footnotesize
\bibliographystyle{apalike}
\bibliography{biblio}

\label{fin}

\end{document}